\documentclass[preprint,12pt]{./elsarticle}

\usepackage{lineno,hyperref}

\usepackage{graphicx}
\usepackage{subcaption}
\usepackage{amssymb}
\usepackage{amsfonts}
\usepackage{amsmath} 
\usepackage{amsbsy} 
\usepackage{amsthm} 
\usepackage{mathrsfs} 
\usepackage{epstopdf} 
\usepackage{dsfont} 

\newtheorem{Theorem}{Theorem}
\newtheorem{hypothesis}{Hypothesis}

\newtheorem{asr}{Condition}

\begin{document}

\begin{frontmatter}

\title{A Shell Bonded to an Elastic Foundation and the Existence of Optimal Elastic and Geometric Conditions for the Membrane Case\tnoteref{t1}}
\tnotetext[t1]{This work is based on a PhD thesis submitted to UCL in 2016 \citep{jayawardana2016mathematical}, and the initial research was funded by The Dunhill Medical Trust [grant number R204/0511] and UCL Impact Studentship.}

\author[mymainaddress]{Kavinda Jayawardana\corref{mycorrespondingauthor}}
\address[mymainaddress]{TEK Optima Research Ltd, Unit 10 Westcroft Buiness Park Oakdene Drive, Three Legged Cross, Wimborne BH21 6FQ}
\ead{kavjayawardana@tekoptimaresearch.com; zcahe58@ucl.ac.uk}
\cortext[mycorrespondingauthor]{Corresponding author}

\begin{abstract}
In this article, we derive a mathematical model for a  shell (i.e. a thin elastic body) bonded to an elastic foundation by modifying Koiter's linear shell equations. We prove the existence and the uniqueness of solutions, and we explicitly derive the governing equations and the boundary conditions for the general case. Finally, with numerical modelling and asymptotic analysis, we show that there exist optimal values for the Young's modulus, the Poisson's ratio and the thickness of the shell (relative to the elastic foundation), and the radius of curvature of the contact region such that the planar solution derived by the shell model (i.e. the membrane case, where stretching effects are dominant) results in a good approximation of the thin body. It is often regarded in the field of stretchable and flexible electronics that the planar solution is mostly accurate when the thin body (i.e. the shell or the membrane) is relatively stiffer (i.e. has a high Young's modulus) than the thicker foundation. As far as we are aware, this is the first analysis showing that for a membrane (or the planar solution of a shell) bonded to an elastic foundation, higher stiffness (relative to the elastic foundation's stiffness) alone would not guarantee a more accurate solution.
\end{abstract}


\end{frontmatter}



\section{Introduction}
Consider a situation where two elastic bodies that are bonded together. For such a case, one can easily model the problem with a simple use of the three-dimensional elasticity equations. Now, consider a scenario where one of the elastic bodies is very thin compared to the other body, and planar in a Euclidean sense (i.e. not curved). Then the thin body can be approximated by a plate or a film. Such models are frequently used in the field of stretchable and flexible electronics. Applications of such models can be found in the field of conformal displays \cite{forrest2004path},  thin film solar cells \cite{choi2008polymers,crawford2005flexible,lewis2006material,pagliaro2008flexible}, electronic skins for robots and humans \cite{wagner2004electronic} and conformable electronic textiles \cite{bonderover2004woven}. For such applications, the extent of the deformation an electronic body can endure, before its basic functions (e.g. conductivity, transparency or light emission) are adversely affected, is immensely important. However, design and process engineers who are working on the implementation of flexible electronics often lack confidence due to a lack of understanding or a lack of input data for reliable modelling tools; thus, there is a tremendous amount of research being conducted in academia as well as in the commercial sector \cite{logothetidis2014handbook}.\\

\begin{figure}[!h]
\centering
\includegraphics[width=1\textwidth]{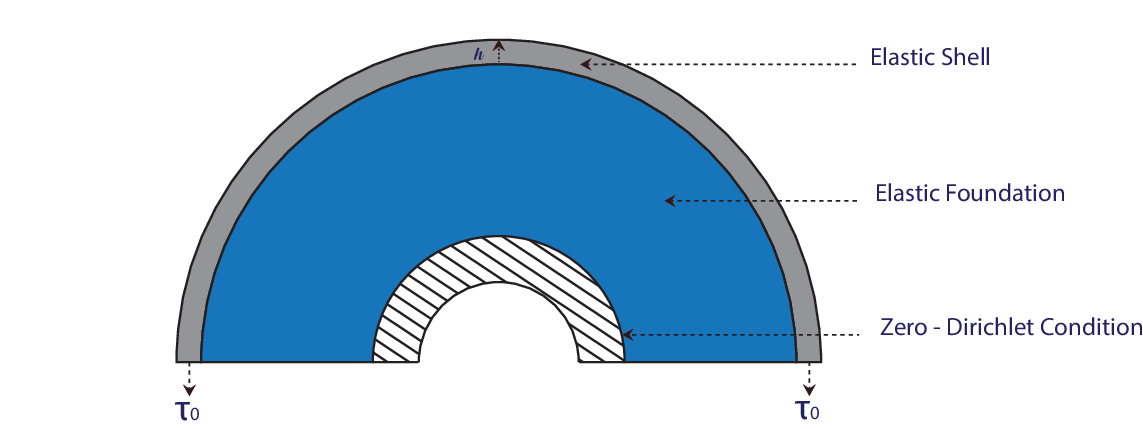}
\caption{Schematic representation of a shell bonded to an elastic foundation, where $h$ is the thickness and $\tau_{\!0}$ is some traction acting on the boundary of the shell}
\label{shellonfoundation}
\end{figure}

Now, suppose that if the intention is to model complex curvilinear stretchable devices for biomedical applications \cite{ko2009curvilinear}, then flat Euclidean models as described above will cease to be useful. Therefore, the goal of our work is to present a simple but a mathematically valid method to model such problems, i.e. a mathematical model for a shell that is bonded to an elastic foundation in a linear-elasticity setting (see Fig. \ref{shellonfoundation}). Also, it is often assumed in the field of stretchable and flexible electronics that the planar solution is mostly accurate when the stiffness of the thin body (i.e. the shell or the membrane) is relatively higher than the stiffness of the thicker foundation \cite{logothetidis2014handbook}, where stiffness means having a high Young's modulus. Baldelli and Bourdin \cite{Andres} attempt to show this asymptotically, and they conclude by stating that the greater the Young's modulus of the thin body relative to the elastic foundation, the more accurate the planar solution. Thus, we also conduct asymptotic and numerical analysis to test this claim.

\subsection{Modelling Difficulties}

A \emph{shell} is a thin 3-dimensional elastic body with a constant thickness $h$ whose rest configuration is curvilinear. Therefore, the governing equation of an elastic shell is derived by modifying the displacement field of this 3-dimensional elastic body by assuming the upper and lower-surfaces of the shell are stress free, and the normal stress is constant through the thickness of the shell, and, finally, integrating the energy functional of the shell from its mid-surface (i.e. $[-\frac12 h, \frac12 h]$ is the interval of integration) and omitting order $h^4$ or above terms to arrive at the governing equations. For a rigorous derivation of the Koiter's shell equations, please consult chapter 4 of Ciarlet \cite{ciarlet2005introduction}. Note that in this article, we mostly adhere to the convention and notations presented by Ciarlet \cite{ciarlet2005introduction}.\\

Now, assume a curvilinear 3-dimensional elastic body with a constant thickness $L$ (where $L$ is much greater than $h$, i.e. $L \gg h$) is bonded to the lower surface of the shell, which we call the \emph{elastic foundation}. Recall that the shell equations are derived by assuming that the both upper and lower-surfaces are stress-free. Therefore, stress present in the contact region (i.e. the existence of stresses in the lower-surface) contradicts the derivation of the shell equations. Also, as the shell equations are integrated through its thickness from its mid-surface, the region $[-\frac12 h, 0)$ will be a area of discontinuity. Now, one could assume the foundation is bonded to the mid-surface of the shell, as in Winkler foundation and Winkler foundation type problems \cite{dillard2018review}. However, this implies that the region $[-\frac12 h, 0)$ is a region of overlap of both the shell and the foundation, resulting in a violation of volume conservation.\\

To overcome such problems, Baldelli and Bourdin \cite{Andres} postulate that one should integrate the energy functional of the thinner body (a plate in authors' work and a shell in our work) from its lower-surface (i.e. $[0, h]$ is now the interval of integration) in their asymptotic study of plates (and films) on elastic Euclidean foundations. The authors further postulate, given that $(\bar\lambda, \bar\mu) \sim  (\varepsilon^q\lambda, \varepsilon^q\mu) $ and $h \sim \varepsilon^p L$, the asymptotic scaling $q-p=2$ results in a model for a film bonded to an elastic  Euclidean foundation and the asymptotic scaling $2<q-p\leq 3$ results in a model for a plate bonded to an elastic Euclidean foundation (i.e. Winkler foundation type model), where $\lambda$ and $\mu$ are the first and  the second Lam\'{e}'s parameters of the thin body respectively, $\bar\lambda$ and $\bar\mu$ are the first and the second Lam\'{e}'s parameters of the elastic Euclidean foundation respectively, $h$ and $L$ are the thicknesses of the thin body and the elastic Euclidean foundation respectively, and $\varepsilon$ is a infinitesimally-small dimensionless parameter (see figure 2 of Baldelli and Bourdin \cite{Andres}). The authors' work implies that stiffer the thin body is relative to the thicker foundation, the more accurate the model. However, Jayawardana \cite{jayawardana2016mathematical} shows that only if one consider the normal and planar displacements as separate cases, then  the asymptotic condition
\begin{align}
\left\{\Lambda h \sim \bar \mu\frac{\ell^2}{L},~ \Lambda h \gg (\bar\lambda+2\bar\mu) L\right\} \label{filmscale}
\end{align}
results in a model for a film bonded to an elastic foundation and the asymptotic scaling
\begin{align*}
\left\{\Lambda h^3\sim \mathscr K \ell^4, ~ \Lambda h^3 \gg L \bar\mu \ell^2 \right\} 
\end{align*}
result in a Winkler foundation type model, where 
\begin{align*}
\Lambda &= 4\mu \frac{(\lambda+\mu)}{(\lambda+2\mu)},
\end{align*}
and where
\begin{align*}
\mathscr K &=  \frac{\bar\mu(3\bar\lambda+2\bar\mu)}{L(\bar\lambda+\bar\mu)}
\end{align*}
is the foundation modulus and $\ell^2$ is the area of the contact region between the thin body and the elastic foundation. Thus, the asymptotic scaling (\ref{filmscale}) implies that for a film bonded to an elastic foundation, higher stiffness alone may not result in a more accurate solution.\\

Thus, for a shell bonded to an elastic curvilinear foundation, following only Baldelli and Bourdin's limits of integration \cite{Andres}, we consider the planar displacements and normal displacement cases separately as follows. For the planar displacements case, i.e. membrane bonded to an elastic curvilinear foundation case, consider a thin overlying body with constant thickness $h$, whose unstrained configuration is described by the diffeomorphism $\bar{\boldsymbol X}\in C^1(\bar\omega\!\times\![0,h);\textbf{E}^3)$, where $\bar{\boldsymbol X} (x^1,x^2,x^3) = \boldsymbol \sigma (x^1,x^2) +x^3 \boldsymbol N (x^1,x^2)$, where $\boldsymbol \sigma \in C^2(\bar\omega ;\textbf{E}^3)$ is an injective immersion that describes the lower-surface of the shell, $\omega \subset \mathbb{R}^2$  is a connected open bounded plane that satisfies the segment condition with a uniform-$C^1(\mathbb{R}^2;\mathbb{R})$ boundary (see definition 4.10 of Adams and Fournier \cite{adams2003sobolev}), $\bar\omega$ defines the closure of $\omega$, $C^k(\cdot)$ is a space of continuous functions that has continuous first $k$ partial derivatives in the underlying domain, $\textbf{E}^k$ is the $k$-dimensional Euclidean space and $\mathbb{R}^k$ is the $k$-dimensional curvilinear space. Now, assume that this overlying body is bonded to an elastic foundation whose unstrained configuration is described by the diffeomorphism $\bar{\boldsymbol X}\in C^1(\omega\!\times\!(-L,0);\textbf{E}^3)$, such that at $x^3 = -L$, the displacement field of the foundation satisfies the zero-Dirichlet boundary condition, i.e. the clamped boundary condition. Thus, in accordance with Baldelli and Bourdin \cite{Andres}, we may approximate the displacement field of both bodies as follows
\begin{align*}
\boldsymbol w = \boldsymbol(w^1,w^2,0\boldsymbol) \mathds{1}_{([0,h])}(x^3) + \left(1+\frac{x^3}{L}\right)\boldsymbol(w^1,w^2,0\boldsymbol)\mathds{1}_{([-L,0))}(x^3),
\end{align*}
where $\boldsymbol(w^1,w^2\boldsymbol) \in  H^1(\omega) \!\times\! H^1(\omega)$ and $\mathds{1}_{(\cdot)}(\cdot)$ is the indicator function (see appendex A.3.vii of Evans \cite{Evans}). Note that unlike Baldelli and Bourdin \cite{Andres}, we make no \emph{prior assumptions} regarding the asymptotic nature of the Young's moduli or Poisson's ratios, or the displacement fields of both bodies. Now, with some asymptotic analysis (without making any prior assumptions), one finds that the condition
\begin{align}
\left\{\Lambda h \sim \bar\mu \frac{\ell^2}{L},~ \Lambda h \gg  (\bar\lambda+2\bar\mu)L, 
~ h B^{\alpha\beta\gamma\delta} F_{\!\text{[II]}\alpha\beta} F_{\!\text{[II]}\gamma\delta}\sim \frac{\bar\lambda+2\bar\mu}{L} \right\} \label{first-scale}
\end{align}
is the only possible asymptotic scaling that yields any valid leading order solutions, where $$\Lambda = 4\mu \frac{(\lambda+\mu)}{(\lambda+2\mu)} ,$$ and where $\lambda$ and $\mu$ are the first and  the second Lam\'{e}'s parameters of the membrane respectively, and $\bar\lambda$ and $\bar\mu$ are the first and the second Lam\'{e}'s parameters of the elastic curvilinear foundation respectively, $\boldsymbol F_{\!\text{[II]}}$ is the second fundamental form tensor of the surface $\boldsymbol\sigma(\omega)$, $\ell^2=\mathrm{meas}(\boldsymbol\sigma(\omega);\textbf{E}^2)$ and $\mathrm{meas}(\cdot;\mathbb{R}^k)$ is the standard Lebesgue measure in $\mathbb{R}^k$ (see chapter 6 of Schilling \cite{schilling2017measures}). Thus, one may collect all the leading order terms (from the energy functional of a shell bonded to an elastic foundation) to find the following energy functional for a membrane bonded to an elastic pseudo-foundation (pseudo-foundation as the displacement of the foundation is now grossly oversimplified) as follows
\begin{align}
J_\text{M}(\boldsymbol w) = ~&\int_{\omega}\left[\frac{1}{2}h B^{\alpha\beta\gamma\delta}\epsilon^M_{\alpha\beta}(\boldsymbol w) \epsilon^M_{\gamma\delta}(\boldsymbol w) + \frac{1}{2}\frac{\bar \mu} {L} w_\alpha w^\alpha - h f^\alpha w_\alpha \right]d\omega\nonumber \\
& - \int_{\partial\omega} h \tau_{\!0}^\alpha w_\alpha ~d(\partial\omega) , \label{Baldelli1}
\end{align}
where $\epsilon^M_{\alpha\beta}(\boldsymbol w)= \frac{1}{2}(\nabla_{\!\alpha} w_\beta + \nabla_{\!\beta} w_\alpha)$ is the strain tensor and $\boldsymbol B$ is the isotropic elasticity tensor of the membrane, $\boldsymbol \tau_{\!0}  \in \boldsymbol L^2(\omega)$ is a traction field (i.e. applied boundary-stress), $\boldsymbol f \in \boldsymbol L^2(\omega)$ is a force density field, $\boldsymbol \nabla$ is the two-dimensional covariant derivative operator with respect to the curvilinear coordinate system $\boldsymbol(x^1,x^2\boldsymbol) \in \omega$, $L^k(\cdot)$ are the standard $L^k$-Lebesgue spaces (see section 5.2.1 of Evans \cite{Evans}). Also, Einstein's summation notation (see section 1.2 of Kay \cite{kay1988schaum}) is assumed throughout, bold symbols signify that we are dealing with vector and tensor fields, and we regard the indices $i,j,k,l \in \{1,2,3\}$ and $\alpha,\beta,\gamma,\delta \in \{1,2\}$. Furthermore, there exists a unique field $\boldsymbol w = \boldsymbol(w^1,w^2\boldsymbol) \in \boldsymbol H^1 (\omega)$ such that $\boldsymbol w$ is the solution to the following minimisation problem
\begin{align*}
J_\text{M}(\boldsymbol w)= \min_{\boldsymbol v \in \boldsymbol H^1 (\omega)} \!\!\! J_\text{M}(\boldsymbol v) ,
\end{align*}
and the unique minimiser $\boldsymbol w$ is also a critical point in $(\boldsymbol H^1 (\omega), J_\text{M}(\cdot))$, where $H^k(\cdot)$ are the standard $W^{k,2}(\cdot)$-Sobolev spaces (see section 5.2.1 of Evans \cite{Evans}). For a rigorous proof and analysis of this problem,  consult section 3.6 of Jayawardana \cite{jayawardana2016mathematical}.\\

As for the normal displacement case, one may again follow Baldelli and Bourdin's approach \cite{Andres} and  define the displacement field as follows
\begin{align*}
\boldsymbol w = \boldsymbol(-x^3\nabla^1 w_3, -x^3\nabla^2w_3, w^3\boldsymbol) \mathds{1}_{([0,h])}(x^3) + \left(1+\frac{x^3}{L}\right)\boldsymbol(0,0,w^3\boldsymbol)\mathds{1}_{([-L,0))}(x^3),
\end{align*}
where $w^3\in  H^2(\omega) $ and naively assume an energy density of the following form
\begin{align*}
W(\boldsymbol w) = ~&\mu \left( \frac{\lambda}{\lambda+2\mu}\rho^\alpha_\alpha(\boldsymbol w)\rho^\gamma_\gamma(\boldsymbol w)  +    \rho^\gamma_\alpha(\boldsymbol w)\rho^\alpha_\gamma(\boldsymbol w) \right) \left(x^3\right)^2 \mathds{1}_{([0,h])}(x^3)\\
& + \frac12\bar\mu \left( \frac{3\bar\lambda+2\bar\mu}{\bar\lambda+\bar\mu}E^3_3(\boldsymbol w)E^3_3(\boldsymbol w) + 4   E^\alpha_3(\boldsymbol w)E^3_\alpha(\boldsymbol w) \right)\mathds{1}_{([-L,0))}(x^3),
\end{align*}
in the hope of deriving a simple asymptotically valid leading order problem, where $\boldsymbol \rho(\cdot) $ is  the change in the second fundamental form tensor of the shell, $\boldsymbol E(\cdot)$ is the strain tensor of the foundation. However, unlike in the Euclidean case, i.e. a plate bonded to an elastic Euclidean foundation case (see section 1.10 of Jayawardana \cite{jayawardana2016mathematical}), now we have $\epsilon_{\alpha\beta}(\boldsymbol w) \mathds{1}_{([0,h])}(x^3) \neq 0$ and $E_{\alpha\beta}( \boldsymbol w) \mathds{1}_{([-L,0))}(x^3) \neq 0$ as a result of the curvilinear nature of our problem. Therefore, one will not be able to express the energy functional as a similar form to equation (\ref{Baldelli1}). Therefore, Baldelli and Bourdin's \cite{Andres} approach fails to work for the normal displacement case for a shell bonded to an elastic foundation.\\

Note that for the mathematical definitions of the second fundamental form tensor, covariant derivative operator, change in the second fundamental form tensor and the strain tensor, we refer the reader to Ciarlet \cite{ciarlet2005introduction}, as it is too cumbersome (i.e. the definitions rely on many other definitions) to define them this early in the text. However, they are mathematically defined in subsequent sections herein in the order they naturally arise.\\

\section{A Shell Bonded to an Elastic Foundation}

Consider an unstrained static three-dimensional isotropic elastic body (which we call the \emph{foundation}) whose volume is described by the diffeomorphism $\bar{\boldsymbol X} \in C^2(\bar\Omega;\textbf{E}^3)$. Now, assume that there exists a thinner isotropic elastic body (which we call the \emph{shell}) bonded to the a subset of the boundary of the foundation such that the contact region is initially stress-free, where this contact region is described by the injection $\boldsymbol \sigma \in C^3(\bar\omega;\textbf{E}^2)$. Thus,  the goal is the find a set of governing equations to describe this problem. Note that in our analysis, we only consider shells that satisfy the following condition:

\begin{asr}\label{dfnShell}
Let the map $\boldsymbol \sigma \in C^3(\bar\omega;\textbf{E}^3)$ describe the contact region between the shell and the foundation, where $\omega \subset\mathbb{R}^2$ is a connected open bounded plane that satisfies the segment condition with a uniform-$C^1(\mathbb{R}^2;\mathbb{R})$ boundary. Given that the thickness of the shell is $h$, we require the following condition to be satisfied
\begin{align*}
0\leq h^2K<  h |H |\ll  1 ~,&~\forall ~ \boldsymbol(x^1,x^2\boldsymbol) \in \bar\omega,
\end{align*}
where $K$ is the Gaussian curvature and $H$ is the mean curvature, i.e. the lower-surface of the shell is not hyperbolic and the thickness of the shell is sufficiently small.
\end{asr}

Prior to our main analysis, we define the following; $T^{ij}(\boldsymbol u) = A^{ijkl}E_{kl}(\boldsymbol u) $ is the second Piola-Kirchhoff stress tensor of the foundation, $E_{ij} (\boldsymbol u) = \frac{1}{2} (\bar\nabla_{\!i} u_j + \bar\nabla_{\!j} u_i)$ is the linearised Green-St Venant strain tensor of the foundation, $$A^{ijkl} = \bar\lambda g^{ij} g^{kl} + \bar\mu (g^{ik}g^{jl} + g^{il}g^{jk})$$ is the isotropic elasticity tensor of the foundation,
\begin{align*}
g_{ij} = \partial_i \bar{\boldsymbol{X}} \cdot \partial_j \bar{\boldsymbol{X}}, ~ i,j \in \{1,2,3\},
\end{align*}
is the covariant metric tensor of $\bar{\boldsymbol X}$,  $\partial_j$ is the partial derivative with respect to the curvilinear coordinate $x^j$,
\begin{align*}
\bar\lambda &= \frac{ \bar\nu \bar E }{(1+\bar\nu)(1-2\bar\nu)}
\end{align*}
and
\begin{align*}
\bar\mu &= \frac{1}{2}\frac{\bar E}{1+\bar\nu}
\end{align*} are the first and the second Lam\'{e}'s parameters of the foundation respectively, $\bar E \in(0,\infty)$ is the Young's modulus of the foundation and $\bar\nu \in(-1,\frac12)$ is the Poisson's ratio of the elastic foundation, $\boldsymbol f$ is an external force density field acting on the elastic foundation and $\bar{\boldsymbol n}$ is the unit outward normal to the boundary $\partial\Omega$ in curvilinear coordinates. Furthermore, $\tau^{\alpha\beta}(\boldsymbol u) = B^{\alpha \beta \gamma \delta} \epsilon_{\gamma \delta}(\boldsymbol u)$ is the stress tensor of the shell, $\eta^{\alpha\beta} (\boldsymbol u) = B^{\alpha \beta \gamma \delta} \rho_{\gamma \delta}(\boldsymbol u)$ is the negative of the change in moments density tensor of the shell, 
\begin{align*}
\epsilon_{\alpha\beta}(\boldsymbol u) = \left[\frac{1}{2} \left(\nabla_{\!\alpha} u_\beta+ \nabla_{\!\beta} u_\alpha\right) - g_3 F_{\!\text{[II]}\alpha \beta} u^3\right]\!|_{\omega}
\end{align*}
is half of the change in the first fundamental form tensor of the shell, $g_3 = \sqrt{g_{33}}|_{\omega} = 1$ by construction, 
\begin{align*}
\rho_{\alpha\beta}(\boldsymbol u) = \Big[g_3\left(\nabla_{\!\alpha} \nabla_{\!\beta} u^3  - F_{\!\text{[II]}\alpha \gamma} F_{\!\text{[II]}\beta}^{~~\gamma} u^3\right) +F_{\!\text{[II]}\beta\gamma}\nabla_{\!\alpha} u^\gamma & \\
+ F_{\!\text{[II]}\alpha\gamma} \nabla_{\!\beta} u^\gamma +\left(\nabla_{\!\alpha}F_{\!\text{[II]}\beta \gamma }\right) u^\gamma & \Big] |_{\omega}
\end{align*}
the change in second fundamental form tensor of the shell,
\begin{align*}
B^{\alpha \beta \gamma \delta} & = \frac{2\lambda\mu}{\lambda + 2\mu}F_{\!\text{[I]}}^{\alpha\beta}F_{\!\text{[I]}}^{\gamma \delta}+ \mu (F_{\!\text{[I]}}^{\alpha \gamma}F_{\!\text{[I]}}^{\beta\delta} + F_{\!\text{[I]}}^{\alpha\delta}F_{\!\text{[I]}}^{\beta\gamma})
\end{align*}
is the isotropic elasticity tensor  of the shell,
\begin{align*}
\lambda &= \frac{ \nu  E }{(1+\nu)(1-2\nu)}
\end{align*}
and
\begin{align*}
\mu &= \frac{1}{2}\frac{ E}{1+\nu}
\end{align*} are the first and the second Lam\'{e}'s parameters  of the shell respectively, $E\in(0,\infty)$ is the Young's modulus  of the shell and $\nu \in(-1,\frac12)$ is the Poisson's ratio of the shell,  $\boldsymbol f_0$ is an external force density field acting on the shell and $\boldsymbol n$ is the unit outward normal vector to the boundary $\partial\omega$ in curvilinear coordinates, and $\boldsymbol \tau_{\!0}$ is an external traction field acting on the boundary of the shell. Finally, $\bar{\boldsymbol \nabla}$ is the covariant derivative operator in the curvilinear space, i.e. for any $\boldsymbol v \in C^1(\bar\Omega;\mathbb{R}^3)$, we define its covariant derivative as follows
\begin{align*}
\bar \nabla_{\!j} v^k = \partial_j v^k + \bar \Gamma^k_{\!ij}v^i,
\end{align*}
where
\begin{align*}
\bar \Gamma^k_{\!ij} = \frac{1}{2}g^{kl}\left ( - \partial_l g_{ij} +\partial_i g_{jl} +\partial_jg_{li} \right)
\end{align*}
are the Christoffel symbols of the second kind, and $\boldsymbol \nabla$ is the covariant derivative operator in the curvilinear plane, i.e. for any $\boldsymbol u \in  C^1(\bar\omega;\mathbb{R}^2)$ we define its covariant derivative as follows
\begin{align*}
\nabla_{\!\beta} u^\gamma = \partial_\beta u^\gamma + \Gamma^\gamma_{\!\alpha\beta}u^\alpha,
\end{align*}
where
\begin{align*}
\Gamma^\gamma_{\!\alpha\beta} = \frac{1}{2}F_{\!\text{[I]}}^{\gamma\delta} \left ( - \partial_\delta F_{\!\text{[I]}\alpha\beta} +\partial_\alpha F_{\!\text{[I]}\beta\delta} +\partial_\beta F_{\!\text{[I]}\delta\alpha} \right)
\end{align*}
are the Christoffel symbols of the second kind in the curvilinear plane.\\

\begin{figure}[!h]
\centering
\includegraphics[width=0.75\textwidth]{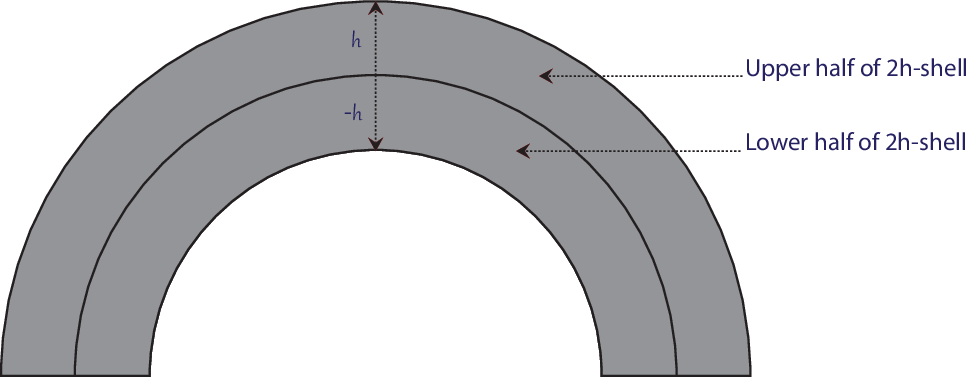}
\caption{Schematic representation of the $2h$-shell}
\label{2h-shell}
\end{figure}

Now, assume for the time being that we are dealing with shell with a thickness $2h$ (see Fig. \ref{2h-shell}) and where, only for the time being, the mid-surface of the shell is described by $\boldsymbol\sigma(\omega)$, and thus, we may express the energy functional of this shell as follows
\begin{align*}
J_{2h}(\boldsymbol u) =~& \int_{\omega} \left[B^{\alpha\beta\gamma\delta}\left(h\epsilon_{\alpha\beta}(\boldsymbol u) \epsilon_{\gamma\delta}(\boldsymbol u) + \frac{1}{3}h^3\rho_{\alpha\beta}(\boldsymbol u) \rho_{\gamma\delta} (\boldsymbol u) \right) - 2hf^iu_i \right] d\omega \\
& - \int_{\partial\omega} 2h \tau_{\!0}^i u_i ~d(\partial\omega),
\end{align*}
where $\boldsymbol u$ describes the displacement field with respect to $\omega$. Note that the energy functional of the $2h$-shell can alternatively be expressed as follows
\begin{align}
J_{2h}(\boldsymbol u)=  2h\bigg (\int_{\omega}\left [ \frac{1}{2}\mu\left(\frac {\lambda }{\lambda+2\mu} \epsilon_{\alpha}^\alpha(\boldsymbol u) \epsilon_{\gamma}^\gamma (\boldsymbol u)  + \epsilon_{\alpha}^\gamma(\boldsymbol u) \epsilon_{\gamma}^\alpha (\boldsymbol u)\right)  -    f^iu_i \right] d\omega & \nonumber\\
-  \int_{\partial\omega}  \tau_{\!0}^i u_i ~d(\partial\omega) & \bigg )\nonumber \\
+ 2h \bigg ( \int_{\omega}\left [\frac{1}{6}\mu\left(\frac {\lambda }{\lambda+2\mu} h^2 \rho_{\alpha}^\alpha(\boldsymbol u) \rho_{\gamma}^\gamma (\boldsymbol u)  + h^2 \rho_{\alpha}^\gamma(\boldsymbol u) \rho_{\gamma}^\alpha (\boldsymbol u) \right) \right] d\omega &\bigg ) .\label{J2h}
\end{align}
Note that as $\boldsymbol \epsilon(\cdot)$ is half of the change in the first fundamental form tensor, we have the following
\begin{align}
\epsilon_{\alpha}^\alpha(\boldsymbol u) \epsilon_{\gamma}^\gamma (\boldsymbol u)   \ll  \frac14 F_{\!\text{[I]}\alpha}^{~\alpha}F_{\!\text{[I]}\gamma}^{~\gamma} &= 1 \label{con1} 
\end{align}
and
\begin{align}
\epsilon_{\alpha}^\gamma(\boldsymbol u) \epsilon_{\gamma}^\alpha (\boldsymbol u) \ll \frac14  F_{\!\text{[I]}\alpha}^{~\gamma}F_{\!\text{[I]}\gamma}^{~\alpha} &= \frac12, \label{con2}
\end{align}
where
\begin{align*}
F_{\!\text{[I]}\alpha\beta} = \partial_\alpha\boldsymbol{\sigma} \cdot  \partial_\beta \boldsymbol{\sigma} ,  ~\alpha,\beta\in \{1,2\},
\end{align*}
is the covariant first fundamental form tensor and $\cdot$ is the Euclidean dot product. Also,  as $\boldsymbol \rho(\cdot)$ is the change in the second fundamental form tensor, condition \ref{dfnShell} implies  the following
\begin{align}
h^2 \rho_{\alpha}^\alpha(\boldsymbol u) \rho_{\gamma}^\gamma (\boldsymbol u)  \ll h^2 F_{\!\text{[II]}\alpha}^{~\alpha}F_{\!\text{[II]}\gamma}^{~\gamma} &= (2hH)^2  \ll 1 \label{con3} 
\end{align}
and
\begin{align}
h^2\rho_{\alpha}^\gamma(\boldsymbol u) \rho_{\gamma}^\alpha (\boldsymbol u) \ll h^2 F_{\!\text{[II]}\alpha}^{~\gamma}F_{\!\text{[II]}\gamma}^{~\alpha} &= 2h^2(2H^2-K)  \ll 1, \label{con4}
\end{align}
where $K = (F_{\!\text{[II]}1}^{~~1}F_{\!\text{[II]}2}^{~~2}-F_{\!\text{[II]}1}^{~~2}F_{\!\text{[II]}2}^{~~1})$, $H=-\frac12F_{\!\text{[II]}\alpha}^{~~\alpha}$,
\begin{align*}
F_{\!\text{[II]}\alpha\beta} = \boldsymbol{N}\cdot \partial_{\alpha\beta}\boldsymbol{\sigma} ,  ~\alpha,\beta\in \{1,2\},
\end{align*}
is the second fundamental form tensor of $\boldsymbol \sigma$,
\begin{align*}
\boldsymbol N = \frac{\partial_1\boldsymbol \sigma \times \partial_2\boldsymbol \sigma}{||\partial_1\boldsymbol \sigma \times \partial_2\boldsymbol \sigma||}
\end{align*}
is the unit normal to the surface $\boldsymbol \sigma$, $\partial_\alpha$ are partial derivatives with respect to curvilinear coordinates $x^\alpha$, $\times$ is the Euclidean cross product and $||\cdot||$ is the Euclidean norm. Now, examining the equations (\ref{con1}) to (\ref{con4}), we may assume that $h^2\rho_{\alpha}^\alpha(\boldsymbol u) \rho_{\gamma}^\gamma (\boldsymbol u) \ll \epsilon_{\alpha}^\alpha(\boldsymbol u) \epsilon_{\gamma}^\gamma (\boldsymbol u)$ and $h^2\rho_{\alpha}^\gamma(\boldsymbol u) \rho_{\gamma}^\alpha (\boldsymbol u) \ll \epsilon_{\alpha}^\gamma(\boldsymbol u) \epsilon_{\gamma}^\alpha (\boldsymbol u)$, and thus, equation (\ref{J2h}) implies that one can expect $J_{2h}(\boldsymbol u)$ to behave approximately linearly in $h$, despite its cubic $h$ dependence, i.e. as result of the geometry and the thinness of the shell, the energy contributions from the bending terms are negligible. This, in turn, implies that the energy stored in the shell's upper and lower halves maybe approximated by dividing the energy functional of the shell by $2$. To be more precise, if $J_{2h}(\boldsymbol u) = J_\text{upper}(\boldsymbol u) + J_\text{lower}(\boldsymbol u)$, then we assume that $\frac{1}{2}J_{2h}(\boldsymbol u) \approx J_\text{upper}(\boldsymbol u) \approx J_\text{lower}(\boldsymbol u)$.   Now, we assume that the  upper half is the form of an overlying shell equation. Thus, we come to the following hypothesis:

\begin{hypothesis} \label{hyp2}
The energy functional of an overlying shell with a thickness $h$ is
\begin{align*}
J_\text{shell}(\boldsymbol u) = ~& \int_{\omega} \left[\frac{1}{2}B^{\alpha\beta\gamma\delta}\left(h\epsilon_{\alpha\beta}(\boldsymbol u) \epsilon_{\gamma\delta}(\boldsymbol u) + \frac{1}{3}h^3\rho_{\alpha\beta}(\boldsymbol u) \rho_{\gamma\delta} (\boldsymbol u) \right) - hf^iu_i \right] d\omega \\
&- \int_{\partial\omega} h \tau_{\!0}^i u_i ~d(\partial\omega),
\end{align*}
given that the shell satisfies the condition \ref{dfnShell}, where $\boldsymbol u$ describes the displacement field with respect to $\omega$.
\end{hypothesis}

To reiterate hypothesis \ref{hyp2} in simple terms, our main postulate is that given a shell with thickness $h$ that is bonded to an elastic foundation, we calculate its energy functional by calculating an energy functional for a shell with a thickness $2h$ and then dividing this energy functional by $2$ (note that this is the same idea used to derive the surface
Cauchy-Bourne model from the Cauchy-Bourne, i.e. dividing the energy contributions from the boundaries of the Cauchy-Bourne model by $2$ to arrive at the surface
Cauchy-Bourne model \cite{jayawardana2013analysis}). Now, with hypothesis \ref{hyp2}, we can express the governing equations for a shell bonded to an elastic foundation as follows:

\begin{Theorem}\label{thrmShell}
Let $\Omega \subset \mathbb{R}^3$ be a connected open bounded domain that satisfies the segment condition with a uniform-$C^1(\mathbb R^3;\mathbb R^2)$ boundary $\partial\Omega$ such that $\omega, \partial\Omega_0\subset\partial\Omega$, where $\bar\omega\cap\bar{\partial\Omega_0} ={\O}$ with $\mathrm{meas}(\partial\Omega_0;\mathbb{R}^2)>0$, and let $\omega \subset \mathbb{R}^2$ be a connected open bounded plane that satisfies the segment condition with a uniform-$C^1(\mathbb{R}^2;\mathbb{R})$ boundary $\partial\omega$. Let $\bar{\boldsymbol X} \in C^2(\bar\Omega;\textbf{E}^3)$ be a diffeomorphism and $\boldsymbol \sigma \in C^3(\bar\omega;\textbf{E}^3)$ be an injective immersion satisfying $0\leq h^2K< h|H|\ll  1$ in $\omega$. Let $\boldsymbol f \in \boldsymbol L^2(\Omega)$, $\boldsymbol f_{\!0} \in \boldsymbol L^2(\omega)$, $\boldsymbol \tau_{\!0} \in \boldsymbol L^2(\partial\omega)$. Then there exists a unique field $\boldsymbol u \in \boldsymbol V_{\!\!\!\mathscr{S}}(\omega,\Omega)$ such that $\boldsymbol u$ is the solution to the following minimisation problem
\begin{align*}
J(\boldsymbol u) = \min_{\boldsymbol v \in \boldsymbol V_{\!\!\!\mathscr{S}}(\omega,\Omega)} \!\!\! \!\!\! J(\boldsymbol v) ,
\end{align*}
where
\begin{align*}
\boldsymbol V_{\!\!\!\mathscr{S}}(\omega,\Omega) = \{ \boldsymbol v \in \boldsymbol H^1(\Omega) \mid ~& \boldsymbol v|_\omega \in H^1(\omega)\!\times\! H^1(\omega)\!\times\! H^2(\omega), \\
~& \boldsymbol v|_{\partial\Omega_0} = \boldsymbol 0, ~ \partial_\beta v^3 |_{\partial\omega} =0 , ~\beta \in\{1,2\}\},
\end{align*} 
\begin{align*}
J(\boldsymbol u) = J_\text{shell}(\boldsymbol u) + \int_{\Omega} \left[\frac{1}{2}A^{ijkl}E_{ij}(\boldsymbol u) E_{kl}(\boldsymbol u) - f^i u_i \right]d\Omega  .
\end{align*}
In particular, the unique minimiser $\boldsymbol u$ is also a critical point in $(\boldsymbol V_{\!\!\!\mathscr{S}}(\omega,\Omega), J(\cdot))$.
\end{Theorem}

\begin{proof}
If the shell is clamped at a subset of its boundary (i.e. $ \boldsymbol u |_{\partial\omega_0} = \boldsymbol 0$ and $ n^\alpha \partial_\alpha u^3 |_{\partial\omega_0} =  0$, where $\partial \omega_0 \subset \partial \omega$, and where this is the displacement-traction problem) then the proof follows trivially from theorem 4.4-3 of Ciarlet \cite{Ciarlet} (and 3.6-1 of Ciarlet \cite{Ciarlet} for the elastic foundation); however, our shell need not be partially-clamped at its boundary, and thus, theorem 4.4-3 of Ciarlet \cite{Ciarlet} fails to be applicable in general. Alternatively, if the shell is not bonded to the elastic foundation (i.e. unique up to a rigid displacement, where this is the pure-traction problem), then the proof follows trivially from theorem 4.4-5 of Ciarlet \cite{Ciarlet}; however, our shell is bonded to the elastic foundation and the foundation is clamped, implying our shell equations are not unique up to a rigid-displacement, and thus, theorem 4.4-5 of Ciarlet \cite{Ciarlet} fails to be applicable. However, notice that we constructed the equations of the bonded shell as a boundary-form (see chapter 4 of Necas \cite{Necas}) to the foundation. Thus, proving the existence and uniqueness of solutions is relatively straight forward as we only need to prove the Korn's inequality by modifying the proof of theorems 3.6-1 and 4.4-3 of Ciarlet \cite{Ciarlet}, which we do as follows.\\

Consider the Hilbert space
\begin{align*}
\boldsymbol W(\omega,\Omega) = \{ \boldsymbol v \in \boldsymbol H^1(\Omega) \mid \boldsymbol v|_\omega \in H^1(\omega)\!\times \!H^1(\omega)\!\times\! H^2(\omega) \},
\end{align*}
equipped with the norm
\begin{align*}
||\boldsymbol v||_{\boldsymbol W(\omega,\Omega) }= \big(||v^1||^2_{H^1(\Omega)} &+\!||v^2||^2_{H^1(\Omega)}\!+\!||v^3||^2_{H^1(\Omega)} \\
&+\!||v^1|_\omega||^2_{H^1(\omega)}\!+\!||v^2|_\omega||^2_{H^1(\omega)}\!+\!||v^3|_\omega||^2_{H^2(\omega)} \big)^{\frac{1}{2}}.
\end{align*}
Now, applying Minkowski inequality (see appendix B.2 of Evans \cite{Evans}) to the Korn’s inequality in curvilinear coordinates without boundary conditions (see theorem 3.8-1 of Ciarlet \cite{ciarlet2005introduction}) and the Korn’s inequality on a surface without boundary conditions (see theorem 4.3-1 of Ciarlet \cite{ciarlet2005introduction}) we obtain the following inequality
\begin{align*}
||\boldsymbol u ||_{\boldsymbol W(\omega,\Omega)} \leq C||\boldsymbol u ||_{\boldsymbol K(\omega,\Omega)},~ \forall~ \boldsymbol u \in \boldsymbol W(\omega,\Omega),
\end{align*}
where
\begin{align*}
\boldsymbol K(\omega,\Omega) = \{ \boldsymbol v \in \boldsymbol L^2(\Omega) \mid ~& \boldsymbol v|_\omega \in L^2(\omega)\!\times \!L^2(\omega)\!\times\! H^1(\omega),\\
& ~\boldsymbol E(\boldsymbol v) \in \boldsymbol L^2(\Omega), ~\boldsymbol \epsilon(\boldsymbol v) \in \boldsymbol L^2(\omega) , ~\boldsymbol \rho(\boldsymbol v) \in \boldsymbol L^2(\omega)\},
\end{align*}
equipped with the norm
\begin{align*}
||\boldsymbol v||_{\boldsymbol K(\omega,\Omega) } = \big (\!||\boldsymbol v||^2_{\boldsymbol L^2(\Omega)}\! &+\! ||v^1||^2_{L^2(\omega)} \!+\!|| v^2||^2_{ L^2(\omega)} \!+\! || v^3||^2_{ H^1(\omega)}  \\
&+\! ||\boldsymbol E(\boldsymbol v)||^2_{\boldsymbol L^2(\Omega)} \! + \!||\boldsymbol \epsilon(\boldsymbol v)||^2_{\boldsymbol L^2(\omega)} \! + \!||\boldsymbol \rho(\boldsymbol v) ||^2_{\boldsymbol L^2(\omega)} \!\big)^{\frac{1}{2}},
\end{align*}
and where $\boldsymbol E(\boldsymbol u) \in \boldsymbol L^2(\Omega)$, $\boldsymbol \epsilon(\boldsymbol u) \in \boldsymbol L^2(\omega) $, $\boldsymbol \rho(\boldsymbol u) \in \boldsymbol L^2(\omega)$ in a sense of distributions, and $C$ is some positive constant that is independent of $\boldsymbol u$.\\

Now, consider the space
\begin{align*}
\boldsymbol V(\omega,\Omega) = \{ \boldsymbol v \in \boldsymbol H^1(\Omega) \mid \boldsymbol v|_\omega \in H^1(\omega)\!\times \!H^1(\omega)\!\times\! H^2(\omega),~\boldsymbol v |_{\partial\Omega_0} =\boldsymbol 0 \},
\end{align*}
equipped  with the norm
\begin{align*}
||\boldsymbol v||_{\boldsymbol S(\omega,\Omega) } = \left(||\boldsymbol E(\boldsymbol v)||^2_{\boldsymbol L^2(\Omega)}\!+\!||\boldsymbol \epsilon(\boldsymbol v)||^2_{\boldsymbol L^2(\omega)}\!+\!||\boldsymbol \rho(\boldsymbol v) ||^2_{\boldsymbol L^2(\omega)} \right)^{\frac{1}{2}},
\end{align*}
where $\boldsymbol E(\boldsymbol u) \in \boldsymbol L^2(\Omega)$, $\boldsymbol \epsilon(\boldsymbol u) \in \boldsymbol L^2(\omega) $ and $\boldsymbol \rho(\boldsymbol u) \in \boldsymbol L^2(\omega)$ in a sense of distributions. The infinitesimal rigid displacement lemma in curvilinear coordinates (see theorem 3.8-2 of Ciarlet \cite{ciarlet2005introduction}) with the condition $\boldsymbol v |_{\partial\Omega_0} =\boldsymbol 0$, and the infinitesimal rigid displacement lemma on a surface (see theorem 4.3-3 of Ciarlet \cite{ciarlet2005introduction}) and the boundary trace embedding theorem (see theorem 5.36 of Adams and Fournier \cite{adams2003sobolev}) with the fact that $\partial \Omega$ is a uniform-$C^1(\mathbb R^3;\mathbb R^2)$ boundary imply that for any $\boldsymbol u \in \boldsymbol V(\omega,\Omega)$, if $\boldsymbol E(\boldsymbol u) = \boldsymbol 0$ in $\Omega$, and $\boldsymbol \epsilon(\boldsymbol u) = \boldsymbol 0 $ and $\boldsymbol \rho(\boldsymbol u) = \boldsymbol 0 $ in $\omega$, then $\boldsymbol u = \boldsymbol 0$ in $\bar \Omega$. Thus, the Rellich-Kondrasov theorem (see theorem 6.3 of Adams and Fournier \cite{adams2003sobolev}) implies that there exists a positive constant $C$ that is independent of $\boldsymbol u$, such that
\begin{align*}
||\boldsymbol u ||_{\boldsymbol W(\omega,\Omega)} \leq C||\boldsymbol u ||_{\boldsymbol S(\omega,\Omega)},~ \forall~ \boldsymbol u \in \boldsymbol V(\omega,\Omega).
\end{align*}

To find weak-solutions, e.g. for finite-element modelling, it is sufficient to conclude the proof here. However, for finite-difference modelling, we require the following final step.\\

Trivial traces (see theorem 5.37 of Adams and Fournier \cite{adams2003sobolev}) and the fact that $\partial \omega$ is a uniform-$C^1(\mathbb{R}^2;\mathbb{R})$ boundary imply that $ \boldsymbol V_{\!\!\!\mathscr{S}}(\omega,\Omega)$ is a proper subset of the space $\boldsymbol V(\omega,\Omega)$ that is closed under the norm $||\cdot ||_{\boldsymbol W(\omega,\Omega)}$. Thus, we can conclude our proof of \emph{Korn's inequality} by confirming the following
\begin{align*}
||\boldsymbol u ||_{\boldsymbol W(\omega,\Omega)} \leq C||\boldsymbol u ||_{\boldsymbol S(\omega,\Omega)},~ \forall~ \boldsymbol u \in \boldsymbol V_{\!\!\!\mathscr{S}} (\omega,\Omega),
\end{align*} 
where $C$ is a positive constant that is independent of $\boldsymbol u$.\\

Now, with the uniform positive-definiteness of the elasticity tensor (see theorem 3.9-1 of Ciarlet \cite{ciarlet2005introduction}) and the uniform positive-definiteness of the elasticity
tensor of the shell (see theorem 4.4-1 of Ciarlet \cite{ciarlet2005introduction}), it is a trivial task to prove that $J(\boldsymbol u)$ is coercive, Fr\'{e}chet differentiable and strictly-convex for $\boldsymbol u$, where they are the conditions imply the existence of unique minimiser that is also a critical point (see theorem 1.5.9 of Badiale and Serra \cite{badiale2010semilinear}). A rigorous proof is found in section 3.4 of Jayawardana \cite{jayawardana2016mathematical}.
\end{proof}

Note that by $ \boldsymbol v|_{\partial\Omega_0}$ and $ \boldsymbol v|_{\omega}$, we mean in a trace sense. Also, the displacement $\boldsymbol u$ restricted to the boundary $\omega$ should be understood in the context of the statement. For example, by $\int_\omega \boldsymbol u~ d\omega$, we mean that $\int_\omega \boldsymbol u|_{\omega} ~d\omega$, where $u|_{\omega}$ is in a trace sense, and we often neglect the term $|_{\omega}$ for convenience.

\subsection{The Equations of Equilibrium}

We assume that $\boldsymbol u \in C^2(\Omega;\mathbb{R}^3)$ with $u^\alpha|_{\omega} \in C^3(\omega)$ and $u^3|_{\omega} \in C^4(\omega)$, and thus, theorem \ref{thrmShell} and the \emph{principle of virtual displacements} (see section 2.2.2 of Reddy \cite{reddy2006theory}) implies that the \emph{governing equations of the elastic foundation} can be expressed as follows
\begin{align*}
\bar\nabla_{\!i} T^i_j(\boldsymbol u) + f_j & =0, ~ j\in\{1,2,3\},
\end{align*}
and the \emph{boundary conditions of the elastic foundation} can be expressed as follows
\begin{align}
\boldsymbol u|_{\partial\Omega_0} &= \boldsymbol 0, \label{ch3clamped}\\
\bar n_i T^i_j(\boldsymbol u)|_{\{\partial\Omega\setminus \{\omega \cup \partial\Omega_0\}\}}&= 0,  ~ j\in\{1,2,3\}, \label{ch3stressfree}
\end{align}
where the Dirichlet boundary condition (\ref{ch3clamped}) is often referred to as the \emph{clamped} or \emph{zero-displacement} boundary condition and the Robin boundary condition (\ref{ch3stressfree}) is often referred to as the \emph{stress-free} boundary condition.\\

Now, the \emph{governing equations of the bonded shell} can be expressed as follows
\begin{align*}
\nabla_{\!\alpha} \tau^\alpha_\beta(\boldsymbol u) &+ \frac{2}{3}h^2F_{\!\text{[II]}\beta}^{~~\alpha}\nabla_{\!\gamma} \eta^\gamma_\alpha(\boldsymbol u) \\
&+ \frac{1}{3}h^2\left(\nabla_{\!\gamma}F_{\!\text{[II]}\beta}^{~~\alpha}\right) \eta^\gamma_\alpha(\boldsymbol u) - \frac{1}{h}\mathrm{Tr}(T^3_\beta(\boldsymbol u)) + f_{0\beta}  = 0 , ~ \beta\in\{1,2\},\\
F_{\!\text{[II]}\alpha}^{~~\gamma}\tau^\alpha_\gamma(\boldsymbol u) &- \frac{1}{3}h^2\nabla_{\!\alpha} \left(\nabla_{\!\gamma}\eta^{\alpha\gamma}(\boldsymbol u)\right) \\
&+ \frac{1}{3}h^2F_{\!\text{[II]}\alpha}^{~~\delta}F_{\!\text{[II]}\gamma}^{~~\alpha} \eta^\gamma_\delta(\boldsymbol u) - \frac{1}{h}\mathrm{Tr}(T^3_3(\boldsymbol u)) + g_3 f_0^3  = 0,
\end{align*}
where $\mathrm{Tr}(T^3_j(\boldsymbol u)) = T^3_j(\boldsymbol u)|_\omega$ and $\mathrm{Tr}(\cdot)$ is the trace operator (see section 5.5 of Evans \cite{Evans}), and the  {boundary conditions of the bonded shell} can be expressed as follows
\begin{align}
\big[n_\alpha\tau^\alpha_\beta(\boldsymbol u) + \frac{2}{3}h^2 n_\gamma F_{\!\text{[II]}\beta}^{~~\alpha} \eta^\gamma_\alpha(\boldsymbol u)\big] |_{\partial \omega} & = \tau_{\!0\beta} ,  ~ \beta\in\{1,2\},\nonumber\\
- \frac{1}{3}h^2 g_3 n _\gamma \nabla_{\!\alpha}\eta^{\alpha\gamma}(\boldsymbol u) |_{\partial \omega} & = \tau_{\!0 3},\nonumber\\
\partial_\beta u^3 |_{\partial\omega} & = 0,~ \beta\in\{1,2\}, \label{Neumann}
\end{align}
where the Neumann boundary condition (\ref{Neumann}) is often referred to as the \emph{zero-slope} boundary condition. Due to the zero-slope boundary condition, boundary moments of the shell is calculated as a part of the solution. Note that on the boundary of the shell (or on a subset of the boundary of the shell, i.e.  on $\partial \omega_0 \subset \partial \omega$), we either have $\partial_\beta v^3 |_{\partial\omega} =0$ (where $n _\alpha\eta^\alpha_\beta(\boldsymbol u) |_{\partial \omega}$ are unknowns) or $n _\alpha\eta^\alpha_\beta(\boldsymbol u) |_{\partial \omega} = \eta_{0 \beta}$ (where $\partial_\beta v^3 |_{\partial\omega}$ are unknowns) for $ \beta\in\{1,2\}$, but never both, where $\boldsymbol \eta_0 \in \boldsymbol L^2(\partial\omega)$ is an external change in moments density field.  If one requires boundary moments of the shell defined prior (i.e. apply external moments as boundary conditions), then one must assume that $\partial_\beta u^3|_{\partial\omega}$ are unknowns. Neglecting the final step of the proof of theorem \ref{thrmShell} implies that a unique solution can be found even if one omits the zero-slope conditions, implying that one can apply boundary moments to the shell, i.e. $n _\alpha\eta^\alpha_\beta(\boldsymbol u) |_{\partial \omega} = \eta_{0 \beta}$ where $\partial_\beta u^3|_{\partial\omega}$ are unknowns. However, from our numerical analysis, we find that omitting the zero-slope boundary conditions leads to an ill-posed  fourth-order finite-difference problem (see section 5.2 of Jayawardana \cite{jayawardana2016mathematical}); thus, we insist upon the zero-slope boundary conditions.

\subsection{A Numerical Example}
\label{ch3Num}

\begin{figure}[!h]
\centering
\includegraphics[width=0.75\textwidth]{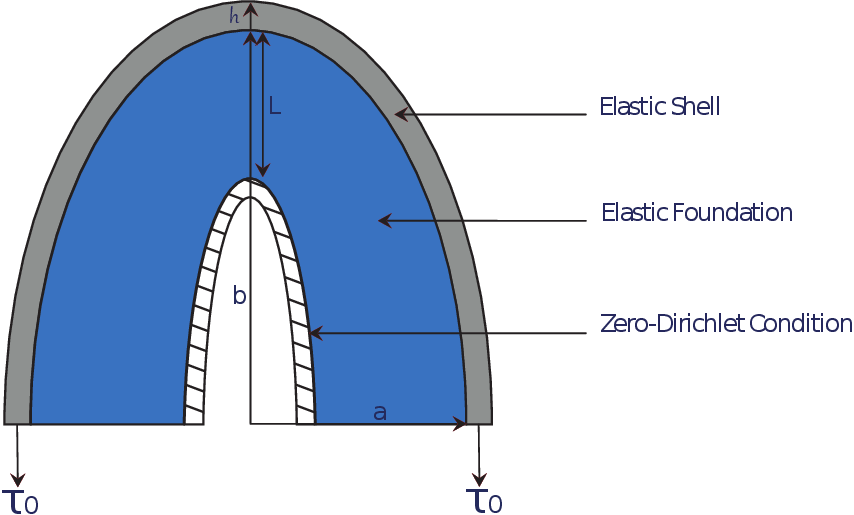}
\caption{Schematic representation of a shell bonded to an elastic foundation, where the cross-section of the contact region forms a semi-ellipse. Note that  $h$ and $L$ are the thicknesses of the shell and the foundation respectively, $b$ and $a$ are the vertical and horizontal radii of the contact region respectively, and $\tau_0$ is the azimuthal traction acting on the boundary of the shell}
\label{shellelliptic}
\end{figure}

To conduct numerical experiments, assume that one is dealing with an overlying shell with a thickness $h$ that is bonded to an elastic foundation, where the unstrained configuration of the foundation defined by an infinitely long annular semi-prism parametrised by the following diffeomorphism $$\bar {\boldsymbol X}(x^1,x^2,x^3) =\boldsymbol(x^1, ~a\sin(x^2), ~b\cos(x^2)\boldsymbol)_\text{E} + \frac{x^3}{\varphi(x^2)}\boldsymbol (0,b\sin(x^2),a\cos(x^2)\boldsymbol )_\text{E} ,$$ where $\varphi (x^2)=(b^2\sin^2(x^2) + a^2\cos^2(x^2))^{\frac{1}{2}}$,  $x^1 \in (-\infty,\infty)$, $ x^2 \in (-\frac{1}{2}\pi,\frac{1}{2}\pi)$, $x^3 \in (-L,0)$, and $a$ is the horizontal radius and $b$ is the vertical radius of the contact region between the shell and the foundation (see Fig. \ref{shellelliptic}). Let $\boldsymbol u = \boldsymbol(0,u^2(x^2,x^3),u^3(x^2,x^3) \boldsymbol )$ be the displacement field, and thus, the governing 
equations of the foundation can be expressed as follows
\begin{align*}
(\bar\lambda + \bar\mu)\partial_2(\bar \nabla_{\!i} u^i ) + \bar \mu (\bar\psi_2)^2 \bar \Delta u^2 & = 0,\\
(\bar \lambda + \bar \mu)\partial_3(\bar \nabla_{\!i} u^i ) + \bar \mu \bar \Delta u^3 & = 0,
\end{align*}
where $\bar\Delta = \bar\nabla_i\bar\nabla^i $ is the vector-Laplacian operator in the curvilinear space (see page 3 of Moon and Spencer \cite{moon2012field}) with respect to $\Omega^\text{New} = (-\frac{1}{2}\pi,\frac{1}{2}\pi) \times (-L,0)$ and where $\bar\psi_2 = \varphi (x^2)+ x^3 ab(\varphi(x^2))^{-2}$.\\

Now, eliminating $x^1$ dependency, the remaining boundaries can be expressed as follows
\begin{align*}
\partial \Omega^\text{New} & = \bar \omega^\text{New}\cup \partial \Omega_0^\text{New} \cup \partial \Omega_f^\text{New},\\
\omega^\text{New} & = (-\frac{1}{2}\pi,\frac{1}{2}\pi)\times \{0\} ,\\
\partial \Omega_0^{\text{New}} & = (-\frac{1}{2}\pi,\frac{1}{2}\pi)\times \{-L\} ,\\
\partial \Omega_f^{\text{New}} &= \{\{-\frac{1}{2}\pi\}\times (-L,0)\}\cup \{\{\frac{1}{2}\pi\}\times (-L,0)\}.
\end{align*}
Thus, the boundary conditions we impose on the foundation reduce to the following
\begin{align*}
u^2|_{\overline{\partial \Omega}_0^\text{New}} &= 0~\text{(zero-Dirichlet)} ,\\
u^3|_{\overline{\partial \Omega}_0^\text{New}} &= 0 ~\text{(zero-Dirichlet)} ,\\
\big[(\bar\psi_2)^2\partial_3u^2 +\partial_2u^3 \big]|_{\partial \Omega_f^\text{New}} & = 0~\text{(zero-Robin)},\\
\big[(\bar\lambda + 2\bar\mu) \partial_2 u^2 + \bar\lambda\left( \partial_3 u^3 + \bar\Gamma^2_{\!22} u^2 + \bar\Gamma^2_{\!23} u^3\right)\big] |_{\partial \Omega_f^\text{New}} & = 0~\text{(zero-Robin)}.
\end{align*}

Now, consider the bonded shell's unstrained configuration that is described by the injective immersion $\boldsymbol \sigma(x^1,x^2) =\boldsymbol (x^1,a\sin(x^2),b\cos(x^2)\boldsymbol )_\text{E}$, where $x^1 \in (-\infty,\infty)$ and $x^2 \in (-\frac{1}{2}\pi,\frac{1}{2}\pi)$. Thus, we may express the governing equations of the shell as follows
\begin{align}
h \Lambda \partial_2 \epsilon^2_2(\boldsymbol u) + \frac{1}{3}h^3 \Lambda (2F_{\!\text{[II]}2}^{~~2}\partial_2 \rho^2_2(\boldsymbol u) + \partial_2F_{\!\text{[II]}2}^{~~2} \rho^2_2(\boldsymbol u)) - \mathrm{Tr} (T^3_2(\boldsymbol u)) & = 0, \label{shellXX1}\\
- h \Lambda F_{\!\text{[II]}2}^{~~2} \epsilon^2_2(\boldsymbol u) + \frac{1}{3}h^3 \Lambda ( \Delta \rho^2_2(\boldsymbol u) - F_{\!\text{[II]}2}^{~~2}F_{\!\text{[II]}2}^{~~2} \rho^2_2(\boldsymbol u)) + \mathrm{Tr} (T^3_3(\boldsymbol u)) & = 0,\label{shellXX2}
\end{align}
where
\begin{align*}
\mathrm{Tr} (T^3_2(\boldsymbol u)) & =\bar\mu\left( (\psi_2)^2 \partial_3u^2 +\partial_2u^3 \right) |_{\omega^\text{New}},\\
\mathrm{Tr} (T^3_3(\boldsymbol u)) & =\big[\bar\lambda\left( \partial_2 u^2 + \bar\Gamma^2_{\!22} u^2 + \bar\Gamma^2_{\!23} u^3\right) + (\bar\lambda + 2\bar\mu) \partial_3 u^3\big]|_{\omega^\text{New}},
\end{align*}
and where $\boldsymbol u |_{\bar\omega^\text{New}} = \boldsymbol (0,u^2(x^2,0),u^3(x^2,0)\boldsymbol )$ is the displacement field of the shell, $\Delta = \nabla_\alpha \nabla^\alpha$ is the vector-Laplacian in the curvilinear plane (see page 3 Moon and Spencer \cite{moon2012field}) with respect to $\omega^\text{New}$ and $ \psi_2 = \varphi(x^2)$.\\

Now, eliminating $x^1$ dependency, we can express the remaining boundaries as follows
\begin{align*}
\partial \omega^\text{New} = \{0\}\cup \{\pi\}
\end{align*}
Thus, the boundary conditions of the shell reduce to the following
\begin{align*}
\big[\Lambda \epsilon^2_2(\boldsymbol u) + \frac23h^2 \Lambda F_{\!\text{[II]}2}^{~~2} \rho^2_2(\boldsymbol u) \big]|_{\{\partial \omega^\text{New}\times [0,h]\}} 
& = \{\tau_0 , \tau_0 \}~\text{(traction)} ,\\
\partial_2 \rho^2_2 (\boldsymbol u)|_{\partial \omega^\text{New} } & = 0~\text{(zero-pressure)} ,\\
\partial_2 u^3 |_{\partial \omega^\text{New} } &= 0~\text{(zero-Neumann)}.
\end{align*}
Despite the fact that the original problem is three-dimensional, it is now a two-dimensional problem as the domain now resides in the set $\{(x^2,x^3) \mid (x^2,x^3)\in [-\frac{1}{2}\pi,\frac{1}{4}\pi]\times[-L,0]\}$.\\

To conduct numerical experiments, we use the second-order accurate fourth-order derivative iterative-Jacobi finite-difference method. Although we use a rectangular grid for
discretisation, as a result of the curvilinear nature of the governing equations, there exists an implicit grid dependence implying that the condition $ \psi_0 \Delta x^2 \leq \Delta x^3$, $\forall ~\psi_0 \in \{\bar\psi_2(x^2,x^3) \mid x^2\in[-\frac{1}{2}\pi,\frac{1}{2}\pi] ~\text{and}~ x^3\in[-L,0] \}$ must be satisfied, where $\Delta x^j$ is a small increment in $x^j$ direction in this context. For our purposes, we let $\Delta x^2 = \frac{1}{N-1}\pi$ and $\psi_0 = \bar\psi_2(\frac{1}{4}\pi,0)$, where $N=250$. We also keep the values $ a = 2$ (units: m), $L=1$ (units: m), $\bar E = 10^3$ (units: Pa), $\bar \nu = \frac{1}{4}$, $\tau_0=1$ (units: N/$\text{m}^2$) fixed for all experiments.\\

\begin{figure}[!h]
\centering
\includegraphics[trim = 2cm 1cm 2cm 1cm, clip = true, width=1\linewidth]{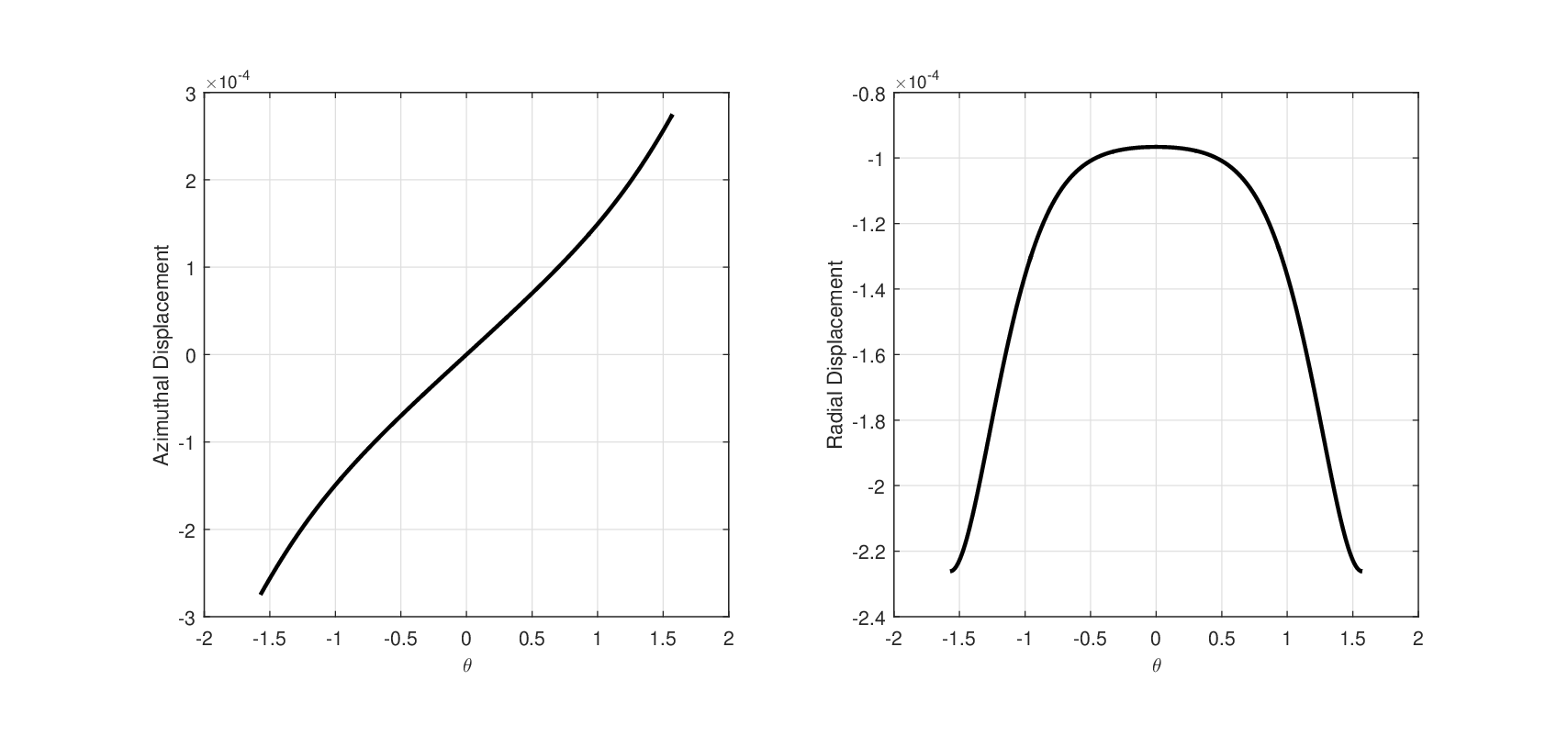}
\caption{Displacement field of a shell bonded to an elastic foundation at its contact region, where $\theta= x^2$\label{Shellch3}}
\end{figure}

Fig. \ref{Shellch3} is calculated with the values $b=2$m, $h=\frac{1}{4}$m, $ E = 6000$Pa and $\nu = \frac{1}{4}$, and it shows the azimuthal (i.e. $u^2$) and the radial (i.e. $u^3$) displacements at the contact region $\omega^\text{New}$. The maximum azimuthal displacements are observed at $x^2=\pm\frac{1}{2}\pi$, with respective azimuthal displacements of $u^2 = \pm 2.75\times 10^{-4}$rad. The maximum radial displacement is observed at $x^2=\pm\frac{1}{2}\pi$, with a radial displacement $u^3 = -2.26 \times 10^{-4}$m. These observations simply imply that the shell (and thus, the foundation) is more likely to deform at the boundaries where we apply external stresses $\tau_0$, and more likely to stay relatively undeformed away from that boundary. Note that all numerical codes are available at \href{http://discovery.ucl.ac.uk/id/eprint/1532145}{http://discovery.ucl.ac.uk/id/eprint/1532145}.

\section{Asymptotic and Numerical Analysis to Find Optimal Elastic and Geometric Properties}

It is asymptotically shown, most notably by Aghalovyan \cite{aghalovyan2015asymptotic}, that a stiffer (i.e. has higher Young's modulus) thin body (relative to the thicker foundation) can result in a more accurate model for a shell bonded to an elastic foundation, where these analyses are conducted by neglecting the planar solution.\\

Furthermore, it is generally regarded in the field of flexible and stretchable electronics that the thin layer must be very stiff relative to the thicker foundation for modelling plates (and films) bonded to elastic foundations (shell or membranes bonded to elastic foundations for our case) to be valid \cite{Andres,logothetidis2014handbook}, where such analyses often concern only the planar effects. However, should one substitute $\boldsymbol(x^2,x^3\boldsymbol) \to \boldsymbol( \bar x^2,L\bar x^3\boldsymbol)$ and $\boldsymbol(u^2,u^3\boldsymbol) \to \boldsymbol(\bar u^2,L\bar u^3\boldsymbol)$ into equations (\ref{shellXX1}) and (\ref{shellXX2}), and examine $\bar u^2$ terms, one finds that the only valid asymptotic scaling of significance is $\phi = 2a \alpha E(e)  \sim 1$, where $E(e)= E(\frac12\pi,e)$ is the complete elliptic integral of the second kind, $E(x^2,e) = \int^{x^2}_0(1-e^2\sin^2(\theta))^\frac12~d\theta$ is the incomplete elliptic integral of the second kind, $e = (1-(b/a)^2)^\frac12$ is the elliptical modulus (see chapter 17 of Abramowitz \emph{et al.} \cite{abramowitz1965handbook}) and $\alpha = (\bar\mu/(hL\Lambda))^\frac12$. This is further justified by our membrane bonded to a foundation model (\ref{Baldelli1}) (which is derived by Baldelli and Bourdin's approach \cite{Andres}), as when applied to the case that we introduced in section \ref{ch3Num}, we find a solution of the following form
\begin{align}
w^2(x^2) &= \frac{\sinh\left(a\alpha E(x^2,e)\right)}{\alpha \Lambda \varphi(x^2) \cosh\left(a \alpha E(e) \right) },\label{AndresSolution}
\end{align}
where $\tau_0 =1$, $\{h \Lambda \sim \bar\mu \frac{\ell^2}{L}, ~  \Lambda h \gg  (\bar\lambda+2\bar\mu)L,~ h B^{22}_{22} F_{\!\text{[II]}2}^{~2} F_{\!\text{[II]}2}^{~2}\sim \frac{\bar\lambda+2\bar\mu}{L}\}$ and\\ $\ell = \mathrm{meas}(\boldsymbol\sigma(\omega^\text{New});\textbf{E})$ (modified as $\omega$ is now no longer bounded). Jayawardana \cite{jayawardana2016mathematical} shows that solution (\ref{AndresSolution}) is only valid when $\phi  \approx 1$ (i.e. when $h \Lambda \sim \bar\mu \frac{\ell^2}{L}$), implying that there exist optimal values for the Young's modulus, the Poisson's ratio and the thickness of the shell (relative to the elastic foundation), and the radius of curvature of the contact region that can result in a more accurate model for a membrane bonded to an elastic foundation (see section 3.6 of Jayawardana \cite{jayawardana2016mathematical}). This contradicts Baldelli and Bourdin’s \cite{Andres} main postulate that the stiffer the membrane (or the film or the plate), then the more accurate the (planar) solution. Thus, our goal in this section is to determine how the relative error between a benchmark model (where the thin layer is not approximated by a shell) and our model for a shell bonded to an elastic foundation behaves for various values of $\delta E =E/\bar E$, $\delta \nu = \nu/\bar \nu$,  $\delta h = h/L$, and $\delta b = b/a$. Note that we assume the default values  $\delta b=1$, $\delta h =\frac18$, $\delta E = 8$ and $ \delta \nu = 1$, unless it strictly says otherwise.\\

For the benchmark model, we numerically model the thin body as a three-dimensional body and we do not approximate this body as a shell. Thus, the displacement field of this bonded two-body elastic problem is obtained by the use of the standard equilibrium equations in the linear elasticity theory.\\

In accordance with the framework that is introduced in section \ref{ch3Num}, the overlying body is restricted to the region $x^3 \in (0,h)$. Thus, we can express the governing equations of the thin body as follows
\begin{align*}
(\lambda + \mu)\partial_2(\bar \nabla_{\!i} v^i ) + \mu (\bar\psi_2)^2\bar \Delta v^2 & = 0,\\
(\lambda + \mu)\partial_3(\bar \nabla_{\!i} v^i ) + \mu \bar \Delta v^3 & = 0,
\end{align*}
where $\boldsymbol v = \boldsymbol(0,v^2(x^2,x^3),v^3(x^2,x^3) \boldsymbol )$ is the displacement field of the thin body, and the following boundary conditions of the thin body
\begin{align*}
\big[(\lambda + 2\mu) \partial_2 v^2 + \lambda\left( \partial_3 v^3 + \bar\Gamma^2_{\!22} v^2 + \bar\Gamma^2_{\!23} v^3\right)\big]|_{\{\{\partial \omega_{T_0}^\text{New}, \partial \omega_{T_\text{max}}^\text{New}\}\times [0,h]\}} 
& = \{\tau_0 , \tau_\text{max}\} ,\\
\big[\lambda\left( \partial_2 v^2 + \bar\Gamma^2_{\!22} v^2 + \bar\Gamma^2_{\!23} v^3\right) + (\lambda + 2\mu) \partial_3 v^3\big]|_{ \{ (-\frac{1}{2}\pi,\frac{1}{2}\pi)\times\{h\}\}} & = 0,\\
\big[(\bar\psi_2)^2\partial_3v^2 +\partial_2v^3\big] |_{\{\partial \omega^\text{New}\times [0,h]\}\cup \{ (-\frac{1}{2}\pi,\frac{1}{2}\pi)\times\{h\}\}} & = 0 ,
\end{align*}
with following equations characterising the bonding of the thin body to the elastic foundation
\begin{align}
\big[u^2-v^2\big]|_{ \omega^\text{New}} &=0 , \label{bc1}\\
\big[u^3-v^3\big]|_{ \omega^\text{New}} & =0 , \label{bc2}\\
\big[\mathrm{Tr} (T^3_2(\boldsymbol u)) - \mu\left((\psi_2)^2\partial_3v^2 +\partial_2v^3\right)\big] |_{ \omega^\text{New}}& = 0 , \label{bc3}\\
\big[\mathrm{Tr} (T^3_3(\boldsymbol u)) - \big( \lambda\left( \partial_2 v^2 + \bar\Gamma^2_{\!22} v^2 + \bar\Gamma^2_{\!23} v^3\right) + (\lambda + 2\mu) \partial_3 v^3\big)\big]|_{ \omega^\text{New}} & = 0  \label{bc4},
\end{align}
where equations (\ref{bc1}) to (\ref{bc4}) represent the continuousness of the azimuthal displacement, the radial displacement, the azimuthal stress and the radial stress respectively.\\

Note that as a result of the grid dependence in the curvilinear coordinates, discretisation of the thin body must satisfy the relation $ \psi_0 \Delta x^2 \leq \Delta x^3$, $\forall ~ \psi_0 \in \{\bar\psi_2(x^2,x^3) \mid x^2\in[-\frac{1}{2}\pi,\frac{1}{2}\pi] ~\text{and}~ x^3\in[0,h] \}$. For our purposes, we let $\Delta x^2 = \frac{1}{N-1}\pi$ and $\psi_0 = \bar\psi_2(\frac{1}{4}\pi,h)$, where $N=250$.\\

\begin{figure}[!h]
\centering
\includegraphics[trim = 2cm 1cm 2cm 1cm, clip = true, width=1\linewidth]{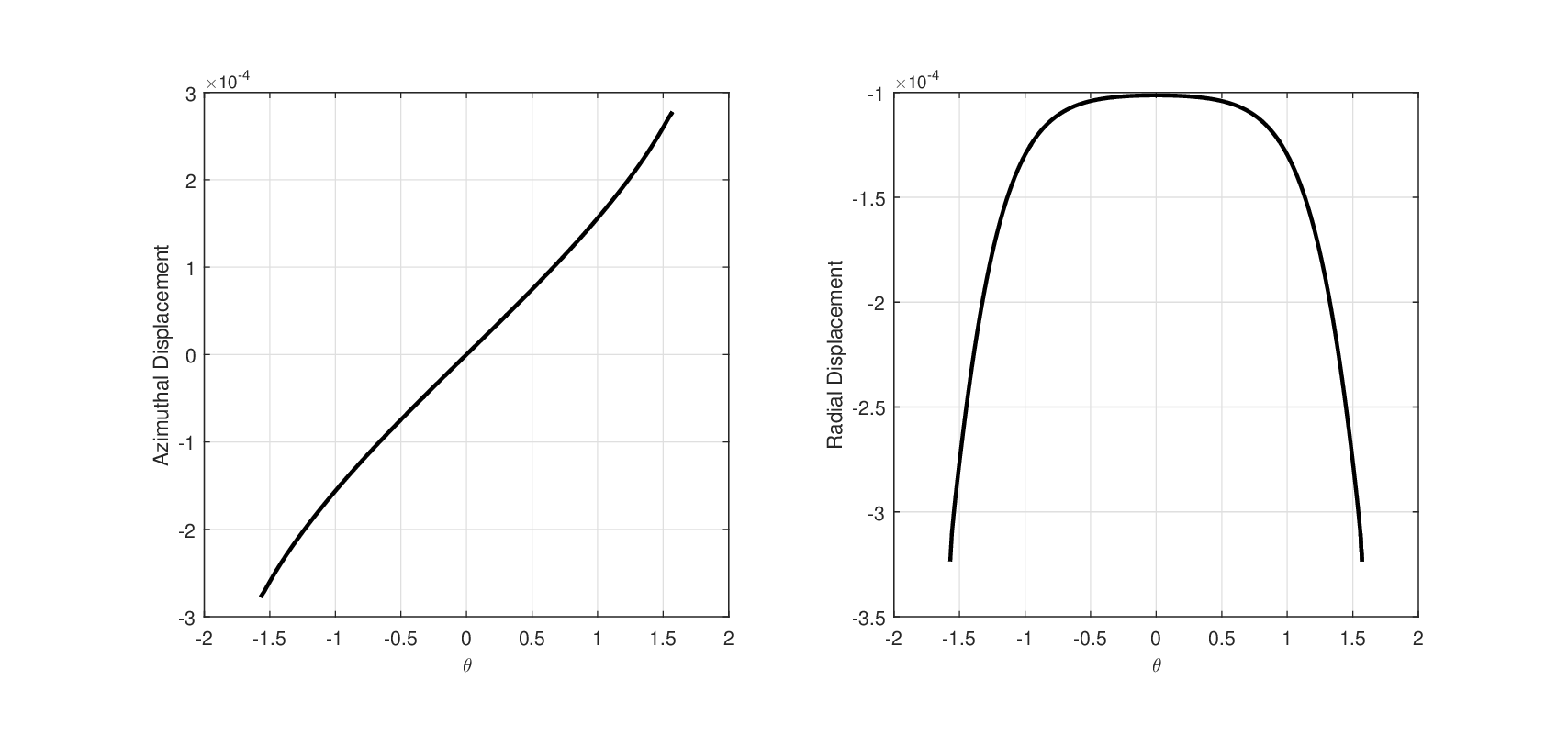}
\caption{Displacement field of the bonded two-body model at the contact region, where $\theta = x^2$\label{Realch3}}
\end{figure}

Fig. \ref{Realch3} is calculated with the values of $b=2$m, $h=\frac{1}{4}$m, $E = 6000$Pa and $\nu = \frac{1}{4}$, and it shows the azimuthal (i.e. $u^2$ and $v^2$) and the radial (i.e. $u^3$ and $v^3$) displacements at the contact region $\omega^\text{New}$. The maximum azimuthal displacements are observed at $x^2=\pm\frac{1}{2}\pi$, with respective azimuthal displacements of $u^2 = v^2 = \pm 2.78\times 10^{-4}$rad. The maximum radial displacement is observed at $x^2=\pm\frac{1}{2}\pi$, with a radial displacement of $u^3 = v^3 = -3.24 \times 10^{-4}$m. Just as it is in the analysis of Fig \ref{Shellch3}, current observations simply imply that the thin body (and thus, the foundation) is more likely to deform at the boundaries where we apply external stresses $\tau_0$, and more likely to stay relatively undeformed away from that boundary.\\

To proceed with our error analysis, we calculate the relative error between the displacement field of the foundation predicted by our model for a shell bonded to an elastic foundation and predicted by the bonded two-body elastic model by the following metric
\begin{align*}
\mathrm{Relative Error}(u^i) &= \frac{\sqrt{\sum_{\{k,l\}}||u_\text{shell-model}^i(y^2_k, y^3_l)-u_\text{two-body}^i(y^2_k, y^3_l)||^2}}{\sqrt{\sum_{\{k,l\}}||u_\text{shell-model}^i(y^2_k, y^3_l)||^2 + ||u_\text{two-body}^i(y^2_k, y^3_l)||^2}}~,
\end{align*}
where $y^2_k = -\frac12\pi +k\Delta x^2$, $y^3_l = -L+l\Delta x^3$, $0\leq k \leq (N-1)$ and   $0\leq l \leq (N-1)$. From this, we should be able to ascertain how stress from the thin body (approximate by a shell or otherwise) propagate to the foundation and deforms it: this, in turn, can help us to identify optimal elastic and geometric conditions that leads to a more accurate model for a shell bonded to an elastic foundation.\\

\begin{figure}[!h]
\centering
\includegraphics[trim = 2cm 1cm 2cm 1cm, clip = true, width=1\linewidth]{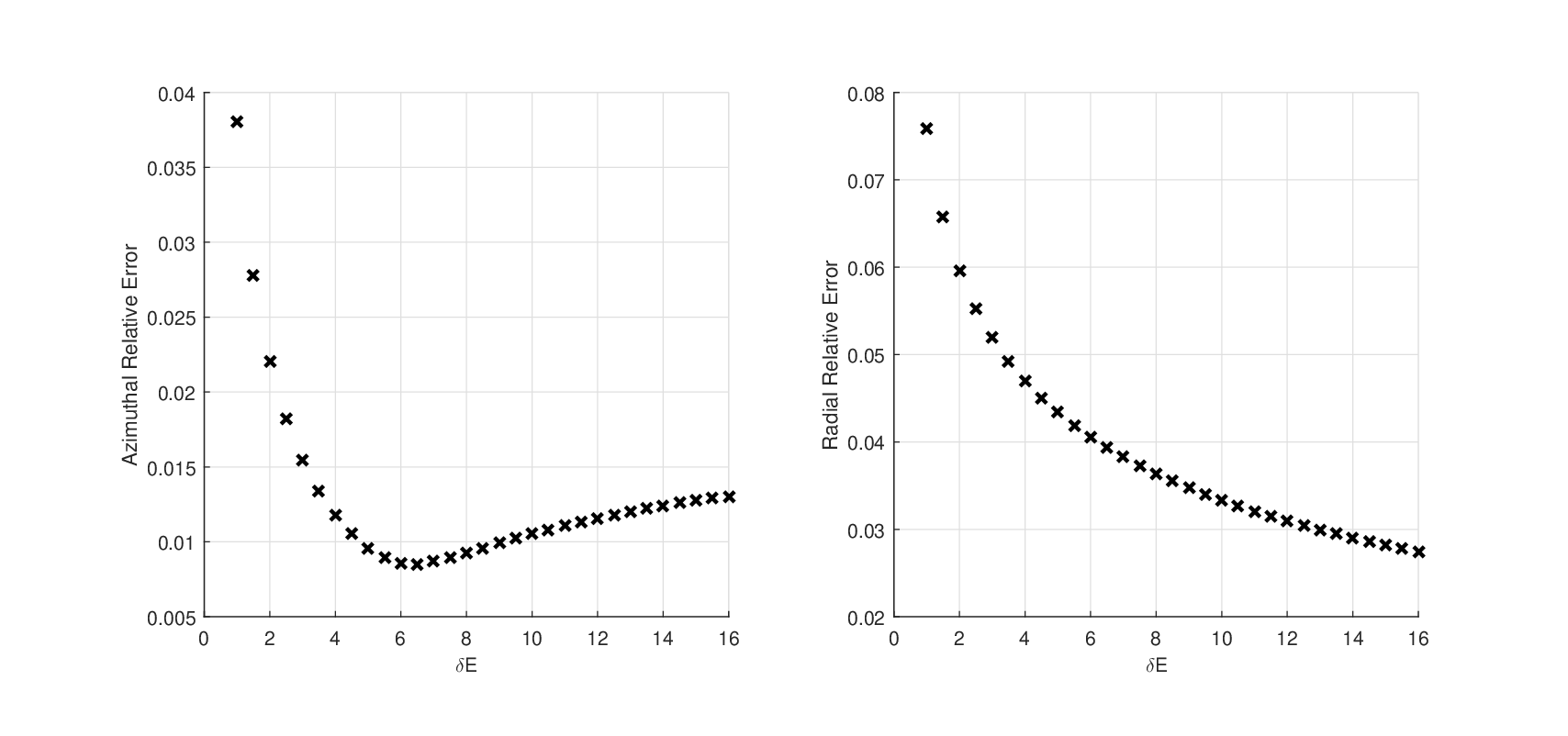}
\caption{Relative error for $\delta E$\label{Ch3EX}}
\end{figure}

\begin{figure}[!h]
\centering
\includegraphics[trim = 2cm 1cm 2cm 1cm, clip = true, width=1\linewidth]{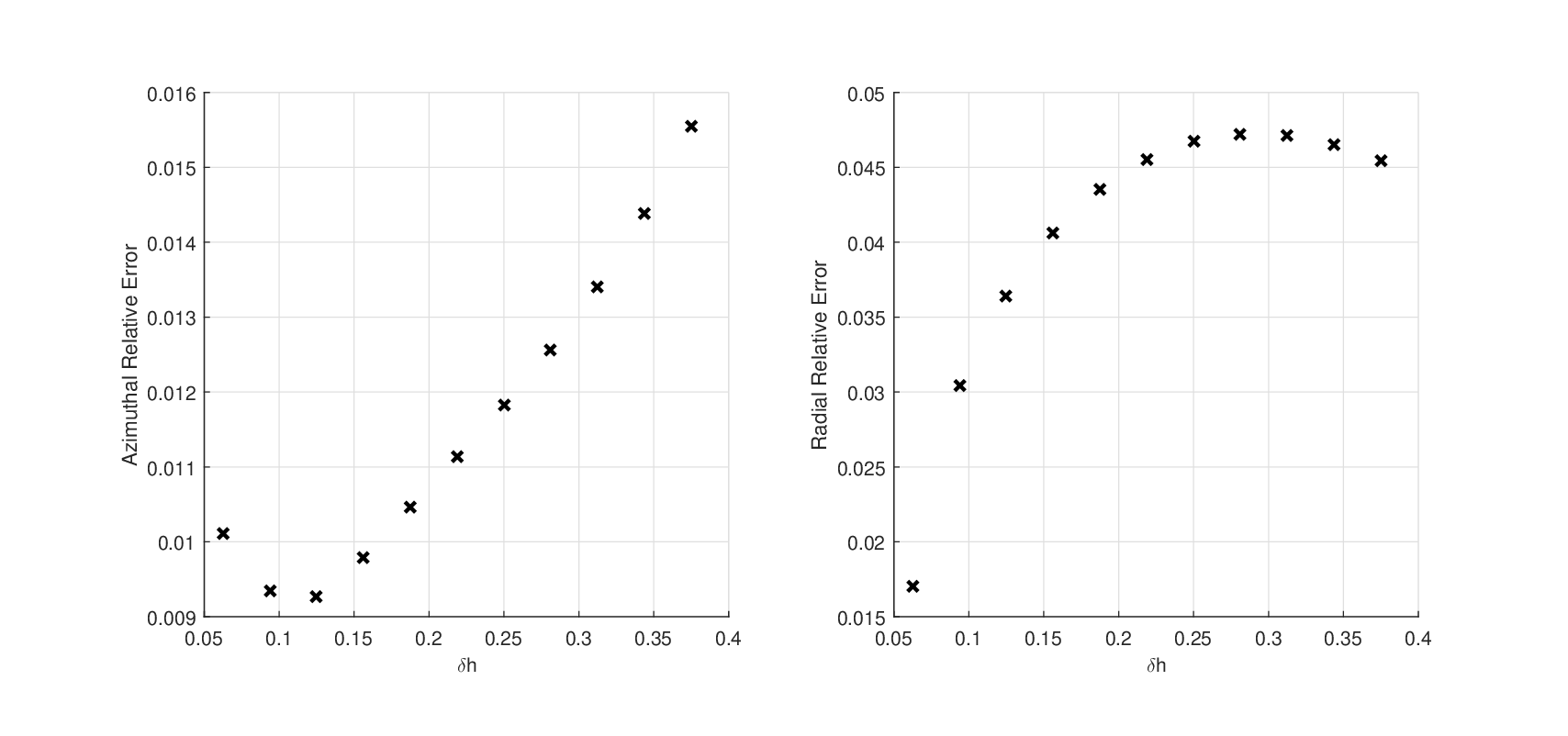}
\caption{Relative error for $\delta h$\label{Ch3hX}}
\end{figure}

\begin{figure}[!h]
\centering
\includegraphics[trim = 2cm 1cm 2cm 1cm, clip = true, width=1\linewidth]{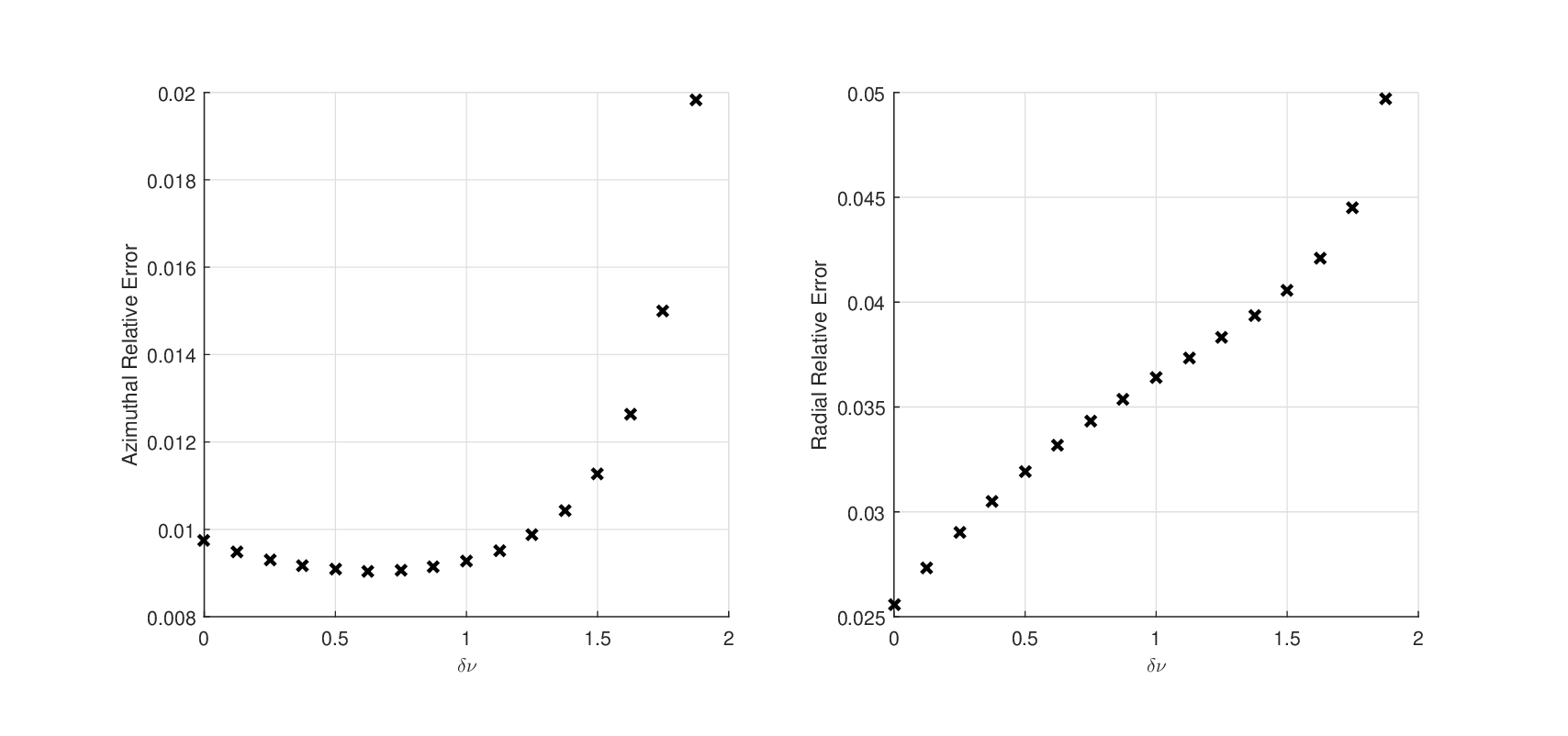}
\caption{Relative error for $\delta \nu$\label{Ch3nuX}}
\end{figure}

\begin{figure}[!h]
\centering
\includegraphics[trim = 2cm 1cm 2cm 1cm, clip = true, width=1\linewidth]{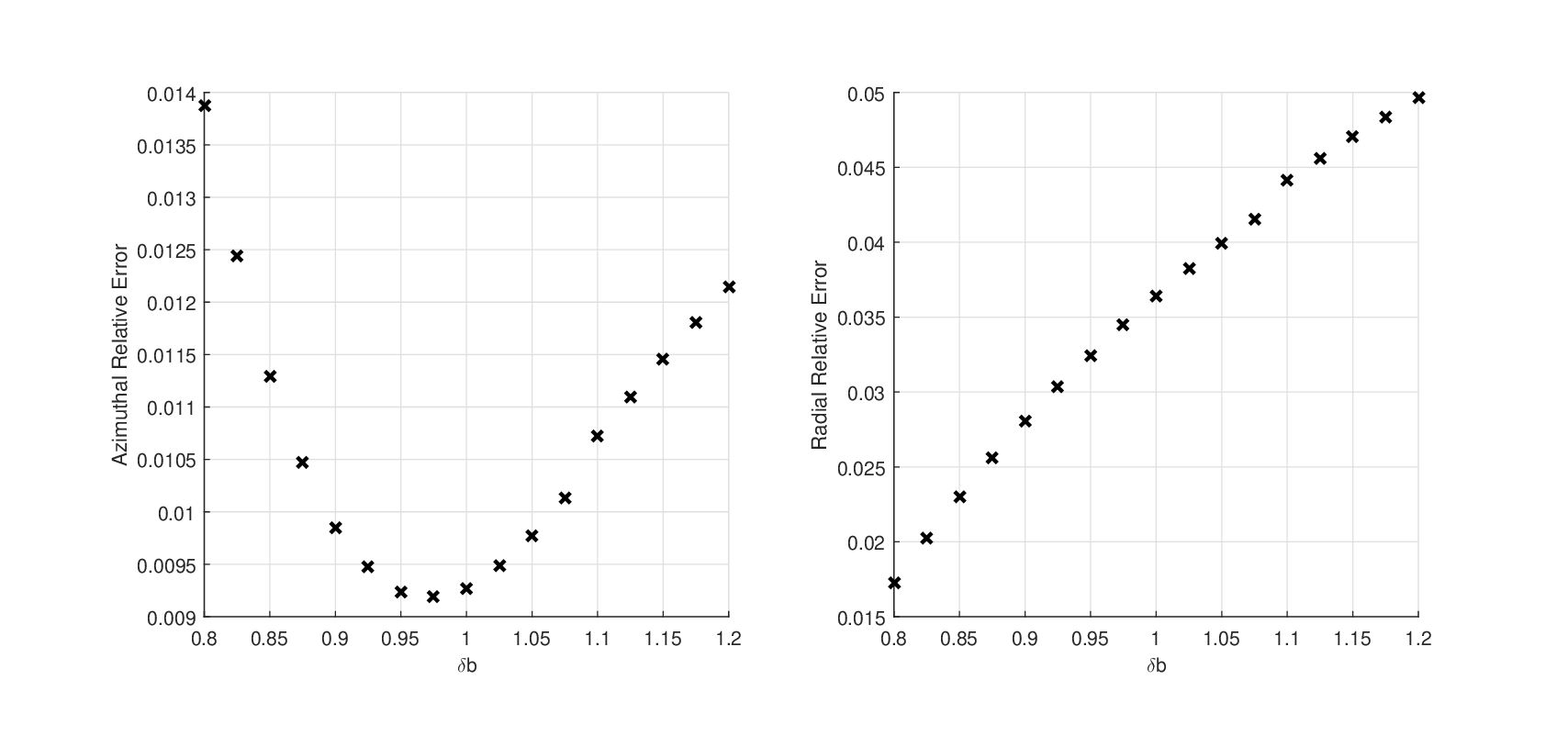}
\caption{Relative error for $\delta b$\label{Ch3cX} }
\end{figure}

Fig. \ref{Ch3EX} shows that the azimuthal error attains a minimum at $\delta E= 6.5$ with a relative error of $0.852\%$. This implies that increasing Young's modulus of the shell indefinitely may not result in the most accurate planar solutions. As for the radial error, one can see that as $\delta E$ of shell increases, the radial error decreases. To reduce the radial error, it appears to be a sound choice to increase the Young's modulus of the shell to relatively high values. Aghalovyan \cite{aghalovyan2015asymptotic} observes similar results for the asymptotic analysis of the modulus of an orthotropic foundation and two-layer anisotropic plates.\\

Fig. \ref{Ch3hX} shows that the azimuthal relative error attains a minimum at $\delta h= 0.125$ with a relative error of $0.927\%$. This implies that there exists an optimum shell thickness where the azimuthal relative error attains a minimum. As for the radial error, one can see that it attains a maximum at $\delta h = 0.25$ with a relative error of $4.67\%$. However, it also shows that the radial error decreases as the thickness of the shell decreases, which is consistent with hypothesis \ref{hyp2}. \\

Fig. \ref{Ch3nuX} shows that the azimuthal error attains a minimum at $\delta \nu= 0.625$ with a relative error of $0.905\%$. This implies that there exists an optimum Poisson's ratio of the shell where the azimuthal relative error is a minimum. As for the radial error, one can see that as $\delta \nu$ of the shell decreases, the relative error also decreases, implying that shell with a relatively low Poisson's ratio may result in lower radial error.\\

Fig. \ref{Ch3cX} shows that the azimuthal error attains a minimum at $\delta b= 0.975$ with a relative error of $0.920\%$. This implies that when the mean curvature of the contact region is almost constant, the azimuthal error attains a minimum. As for the relative radial error, one can see that as $\delta b$ decreases, the relative error also decreases. Assuming a contact region $x^2\in(-\frac{1}{2}\pi,\frac{1}{2}\pi)$ (with $x^2\neq\pm\frac{1}{2}\pi)$, the latter observation may be interpreted as follows: as the radius of curvature of the contact region increases, the radial error decreases. This appears to be consistent with hypothesis \ref{hyp2}, as we derived our bonded shell equation to be valid for contact regions with a high radius of curvature (see condition \ref{dfnShell}).\\

Fig.s \ref{Ch3EX} to \ref{Ch3cX} imply the existence of optimal conditions where azimuthal error attains a minimum, and all four cases coincide with our asymptotic scaling $\phi  \sim 1$ as we observe minimum errors for elastic and geometric values implied by the condition $\phi  \approx 1$. Given that the elastic foundation has a constant thickness $L$ and the lower boundary of the foundation satisfies the zero-Dirichlet boundary condition, this asymptotic scaling can be expressed as follows $$ \Lambda h \sim \bar\mu \frac{\ell^2}{L} ,$$ which is implied by asymptotic condition (\ref{first-scale}).

\section{Conclusions}

In our analysis, we studied a shell bonded to an elastic foundation. We derived our mathematical model by modifying Koiter's linear shell equations (see Ciarlet \cite{Ciarlet}). Then, we used Ciarlet's  work  \cite{ciarlet2005introduction} to prove the existence and the uniqueness of solutions, and we explicitly derived the governing equations and the boundary conditions for the general case for a shell bonded to an elastic foundation. Although, we have shown the existence and the uniqueness of weak solutions, we did not prove any higher regularity results, which are vital for convergence of numerical solutions. This we leave for future work.\\

For numerical analysis, we conducted error analysis to see how well our model for a shell bonded to an elastic foundation can approximate the displacement field of the foundation of a two-body contact problem modelled with the standard equilibrium equations in linear elasticity. Our analysis shows that the radial solution of our bonded shell model can approximate the displacement field of the foundation with a significant degree of accuracy given that the Young's modulus of the shell is high, which is consistent with analogous models that exist in the literature \cite{aghalovyan2015asymptotic}. However, both our numerical and asymptotic analyses show that there exist optimal values for the Young's modulus, the Poisson's ratio and the thickness of the shell (relative to the foundation), and the radius of curvature of the contact region where we observe a minimum azimuthal error. Our numerical modelling also implies that the radial error is a minimum for a shell if it has a relatively low Poisson's ratio and is relatively thin, and if the contact region (between the shell and the elastic foundation) has a high radius of curvature, where the latter two conditions are consistent with the derivation of our model. \\

It is often regarded in the field of stretchable and flexible electronics that the planar solution (where stretching effects are dominant) is mostly accurate when the stiffness of the thinner body (e.g. a plate, a shell, a film or a membrane) is relatively higher than the thicker elastic foundation (i.e. the greater the Young's modulus of the thin body relative to the elastic foundation, the more accurate the planar solution) \cite{Andres,logothetidis2014handbook}. The significance of our work is that, as far as we are aware, this is the first analysis conducted on the planar solution  (both asymptotically and numerically) showing that indefinitely increasing the stiffness of the thinner body will not guarantee a more accurate solution, as there exists an optimum Young's modulus (also other optimal elastic and geometric properties) that can result in a more accurate model for a membrane (or a film) bonded to an elastic foundation.\\

On a final note, the motivation to numerically model the azimuthal and and radial error separately is to show that the elastic and geometric conditions that result in an accurate planar solution (a shell-membrane or a film on a elastic foundation) is vastly different to conditions result in an accurate normal solution (only the radial component of a shell or a plate on an elastic foundation). Thus, assuming the optimal conditions that is valid for the normal solution (e.g. large Young's modulus) is also applicable to the planar solution (i.e. the assumption common in the field of stretchable and flexible electronics), may not be a justifiable assumption.\\

\section*{Acknowledgments}
We thank Dr Nick Ovenden (UCL) and Prof Alan Cottenden (UCL) for their supervision, Brad Turner (TEKOR) for his assistance, and Christopher Law for the illustrations.

\bibliographystyle{./model1-num-names}
\bibliography{ShellsOnElasticFoundations}%
\biboptions{sort&compress}

\end{document}